\documentclass[11pt]{article}
\catcode`\@=11 \catcode`\@=12
\usepackage[utf8]{inputenc}
\usepackage[T1]{fontenc}
\usepackage{amsmath}
\usepackage{amssymb}
\usepackage{amsthm}
\usepackage{oldgerm}
\usepackage{multicol}
\usepackage{verbatim}
\usepackage{color}
\usepackage{ulem}
\usepackage{enumitem}
\usepackage{mathrsfs}
\usepackage{graphicx}
\usepackage{dutchcal}

\usepackage[pdftex,colorlinks,urlcolor=blue,pdfstartview=FitH,
]{hyperref}

\makeatletter
\newcommand{\mylabel}[2]{#2\def\@currentlabel{#2}\label{#1}}
\makeatother

\newcommand{\Rm}{\mathbb{R}}

\newcommand{\mL}{\mathcal{L}}

\newcommand{\gS}{\ensuremath{\mathfrak{S}}}
\newcommand{\gs}{\ensuremath{\mathfrak{s}}}
\newcommand{\gR}{\ensuremath{\mathfrak{R}}}

\newcommand{\gL}{\ensuremath{\mathfrak{L}}}
\newcommand{\ms}{\ensuremath{\mathcal{s}}}
\newcommand{\mC}{\ensuremath{\mathcal{C}}}

\newcommand{\mP}{\ensuremath{\mathcal{P}}}

\newcommand{\Zm}{\ensuremath{\mathbb{Z}}}

\newcommand{\mM}{\ensuremath{\mathcal{M}}}

\newcommand{\mW}{\ensuremath{\mathcal{W}}}

\newcommand{\mE}{\ensuremath{\mathcal{E}}}

\newcommand{\vs}{\vspace{.2cm}}
\newcommand{\Sp}{\textnormal{\textbf{Sp}}}
\newcommand{\A}{\textnormal{\textbf{A}}}
\newcommand{\M}{\textnormal{\textbf{M}}}

\newtheorem{lem}{Lemma}

\newtheorem{thm}{Theorem}

\newtheorem{cor}[lem]{Corollary}
\newtheorem{prop}[lem]{Proposition}

\newtheorem{defn}[lem]{Definition}

\def\proof { \vspace{.2cm} \noindent{\sc{Proof. }}}
\def\qed {\mbox{}\hfill {\small \fbox{}} \\}
\def\lto{\longrightarrow}
\def\lmto{\longmapsto}

\def\leq{\leqslant}
\def\geq{\geqslant}

\date{}

\begin{document}

\begin{center}
	\begin{huge}
		{\bf Non-degeneracy of closed orbits for generic potentials.}\\
	\end{huge}
	\vs \vs \vs

\noindent
Patrick Bernard,\\
PSL Research University,\\
Université Paris-Dauphine,\\
CEREMADE (UMR CNRS 7534)\\
75775 PARIS CEDEX 16,  
France\\
\texttt{patrick.bernard@dauphine.fr}\\

\end{center}

\vs

%
%
%
%
%
%
%
%
%
%
%
%
%
%
%
\vspace{1cm}
\textit{Résumé :}
On démontre qu'un Hamiltonien convexe générique au sens de Mañé n'a que des orbites périodiques non-dégénérées sur un niveau d'énergie donné. Ce résultat a déjà été énoncé, mais pas démontré, dans la littérature.

\textit{Abstract :}
We prove that Mañé generic convex Hamiltonians have only non-degenerate periodic orbits on a given energy level. 
This result was stated, but not proved, in the literature.

\vspace{1cm}

It is known that a generic dynamical system has only hyperbolic periodic orbits, this is part of the Kupka-Smale theorem, see 
\cite{K63,S63}. When some additional structure is preserved the situation can be different. Here we shall focus on Hamiltonian systems.
Such systems can stably have non-hyperbolic periodic orbits.
This can be understood at  the linear level: hyperbolic matrices are not dense in symplectic matrices.
Yet the periodic orbits of generic Hamiltonian systems satisfy some properties, in particular they are non degenerate (on a given energy level), see \cite{CR1,CR2}. The Bumpy metric theorem states that this is also true for generic geodesic flows, see \cite{A68,A82}.
Mechanical systems, defined on a cotangent bundle from a Hamiltonian
which is the sum of a kinetic energy (a Riemannian metric) and  a potential energy (a function on the base)
also form a natural case of study.
For such systems, one can ask what properties hold for a generic potential once the metric has been fixed
(this is the genericity in the sense of Mañé).
One of the outcomes of the present paper is that the answer  is actually the same as for general Hamiltonian : for a generic potential, all periodic orbits on a given energy surface are non-degenerate. 
This  is not surprising, and actually follows from results previously stated in the literature, see \cite{O08,RR}.
However, the existing proofs apply  only to energy levels higher than the supremum of the potential energy but not below.
The proof in \cite{O08} follows rather closely Anosov's proof of the Bumpy Metric Theorem
with one different step, the perturbation Lemma describing to what extent the linearized flow along a given orbit can be modified by a perturbation preserving this orbit. This perturbation Lemma is obtained  in \cite{O08} in dimension 2, and then in \cite{RR}
(with a correction in \cite{AB}) in higher dimension, elaborating on
a method introduced in \cite{Co10} in the context of geodesic flows. 

However, it is well known that the dynamics of a mechanical system on energy levels lower than the supremum of the potential is 
significantly different  from the one of a geodesic flow, and its study requires different methods. As to the question we are considering, the difficulty comes from the possible presence of what we will call round trip orbits (also sometimes called librations)
following back and forth a same path on the base manifold. 
These orbits fall outside of the scope of the previous papers,
and the goal of the present work is to study them.
The perturbation Lemma for these orbits is harder than for geodesic-like orbits because 
a potential will necessary perturb the two branches of the orbit, and one has to understand how these two perturbations interact.
In the case of mechanical system, the reversible structure allows to tame these interactions, and we manage to prove a perturbation Lemma in that setting.

 Round trip orbits can also exist in the more general class of 
Hamiltonian systems convex in the fibers. In general there is no reversible structure associated to the linearized flow along such an orbit, and the perturbation Lemma remains open.
 To solve this difficulty, we prove that a round trip orbit either admits some kind of reversibility which allows to prove a perturbation Lemma, or can be turned to a geodesic-like orbit  by the addition of a small potential.
This allows us to provide a full proof of the result previously stated, but unproved, in the literature claiming that a Mañé generic convex Hamiltonian has only nondegenerate periodic orbits on a given energy level.

\section{Introduction}

We study the dynamics of  convex Hamiltonian systems on the cotangent bundle of a smooth manifold $M$ of dimension $d+1$. 
In the present paper, smooth means $C^{\infty}$ and   convex Hamiltonian  means that the fiberwise Hessian 
of the Hamiltonian function is positive definite at each point. 
We usually denote by $H$ the Hamiltonian, by $V_H$ the corresponding Hamiltonian vector-field and by $\varphi$ the associated flow. We will focus on the generic properties, in the sense of Ma\~né, of the restricted linearized return map of periodic orbits.
 More precisely, we can take a Poincaré section
of the flow restricted to the energy level of the orbit, and consider the differential at the orbit of the return map, this is what we call the restricted linearized return map.
It is, up to a choice of coordinates, an element of the symplectic group $\Sp(2d)$ if
the manifold $M$ has dimension $d+1$. 
The restricted linearized return maps associated to different sections or to different coordinates are conjugated in $\Sp(2d)$, so the conjugacy class of the restricted linearized return map of a given orbit is well-defined.
Let us define $\Upsilon \subset \Sp(2d)$ as the space of matrices having an eigenvalue equal to a root of the unity or having a double eigenvalue. This is a  countable union of closed sets with empty interior  in $\Sp(2d)$, invariant under conjugacy.
Thanks to this conjugacy invariance, the property, for a periodic orbit, to have a restricted linearized return map in 
$\Upsilon$ or out of $\Upsilon$ is intrinsically meaningful. 
Our main result is :

\begin{thm}\label{thm-g1}
	Given a smooth ($C^{\infty}$) convex Hamiltonian $H: T^*M\rightarrow \Rm$, 
	the following property is satisfied for a generic potential  $u\in C^{\infty}(M)$ :

	The zero energy level of $H+u$ is regular, and it does not contain any periodic orbit having its restricted linearized return map in $\Upsilon$.
\end{thm}

As usual, generic means that the set of potentials which do not satisfy the conclusion is  a  countable union of nowhere dense sets. As is well known, $C^{\infty}(M)$ with its usual topology is a separable Fréchet space, 
meaning that its topology can be defined by a translation invariant complete distance with convex balls, and that it is separable for this distance. In particular it satisfies the Baire property: a countable union of nowhere dense sets has empty interior.

Theorem \ref{thm-g1} is stated in \cite{O08} in dimension $2$, and extended in  \cite{RR}  to all dimensions, see also  \cite{Co10}
where some important ideas are introduced. However, as was already noticed in \cite{AB2}, 
there is a missing case in these works. 
In order to explain the difficulty, we introduce the following definition: 

We say that the periodic orbit $\theta(t)=(Q(t), P(t))$ of minimal period $T$ is neat if there exists a non-trivial open interval $I\subset \Rm/T\Zm$ of times  such that $Q$ is an embedding of $I$ into $M$ and such that 
$Q(I)$ does not intersect $Q(\Rm/T\Zm-I)$.

The methods of \cite{O08,RR} apply only to neat orbits, see Section \ref{sec-rmr}. So the statement actually following from the proofs given in  \cite{O08}  and \cite{RR}
(completed in \cite{AB}) is :

\begin{thm}\label{thm-g2}
Given a smooth convex Hamiltonian $H: T^*M\rightarrow \Rm$, 
the following property is satisfied for generic $u\in C^{\infty}(M)$ :

The zero energy level of $H+u$ is regular, and all \textbf{neat} 
 periodic orbits of $H+u$ which have zero energy  have their restricted linearized 
return map in the complement of $\Upsilon$.
\end{thm}

A formal proof of the statement in this form is given in \cite{AB2} where the non convex case is also discussed. Any
conjugacy invariant set  $\Upsilon\subset \Sp(2d)$ which is a countable union of nowhere dense sets can be considered here instead of the explicit one defined above.

However, not all periodic orbits are neat.
In the case of a convex Hamiltonian (which is the only case we consider in the present paper) the orbits which are not neat are very specific :

\begin{prop}\label{prop-neat}
	Let $H$ be a smooth convex Hamiltonian, and let $\theta(t)=(Q(t),P(t))$ be a periodic orbit of minimal period $T$.
	
	If $\theta$ is not neat, then it is a \textit{round trip} orbit, meaning that 
	there exists an orientation reversing diffeomorphism $\sigma$ of $\Rm/T\Zm$ such that $Q(\sigma(t))=Q(t)$ for each $t$. This diffeomorphism then has two fixed points $\nu_0$ and $\nu_1$, and $Q'(\nu_i)=0$.
	\end{prop}

Below, we will also call $\sigma$ a lifting of $\sigma$ to a decreasing diffeomorphism of $\Rm$.
This Proposition is proved  in section \ref{sec-proj}.
In the present paper, we prove the missing part of Theorem \ref{thm-g1}, namely :

\begin{thm}\label{thm-g3}
	Given a smooth  ($C^{\infty}$) convex Hamiltonian $H: T^*M\rightarrow \Rm$, 
	the following property is satisfied for generic $u\in C^{\infty}(M)$ :

All \textbf{round trip} 
	 orbits of $H+u$ which have zero energy  have their restricted linearized 
	return map in the complement of $\Upsilon$.
\end{thm} 

It is not completely clear to what extent this statement can be extended to different subsets  $\Upsilon \subset \Sp(2d)$. Clarifying this question would require a finer study of reversible symplectic matrices than the basic one provided in  Proposition \ref{prop-gen} below.

The first intuition concerning round trip orbits might be that they do not usually exist or more precisely that $H+u$ does not admit
any round trip periodic orbit for generic $u$. However, this depends on $H$. For example, if $H$ is reversible, meaning that 
$H(q,-p)\equiv H(q,p)$, then round trip orbits, sometimes called librations in that setting, do often exist. 
They are studied for example in \cite{BK,K76}, and the results in these papers imply that librations exist for open sets of potentials, so their study can't be avoided.

We denote by $\check M$ the set of points $q\in M$ such that the restriction of $H$ to $T^*_qM$ has a minimum,
and then we denote by $(q, \wp(q))$ this minimum. The fiberwise convexity (with positive definite Hessian) implies  
that $\check M$ is open in $M$, and that $\wp$ is a smooth section of $T^*\check M$.
 We denote by $\Gamma\subset T^*M$ its graph. In the reversible case, $\Gamma$ is the zero section.
 All round trip orbits are contained in $T^*\breve M$, so there would not be much loss of generality in assuming that $\breve M=M$.

Given $p\in T^*_q\breve M$, $p\neq \wp(q)$, there is one and only one point 
$\tilde p\neq p\in  T^*_qM$ such that $\partial_pH(q,p)$ and $\partial_pH(q, \tilde p)$ are proportional and $H(q,\tilde p)=H(q,p)=:e$. This point is the minimizer of the linear form
$\partial_p H(q,p)$ on $T^*_qM\cap \{H= e\}$. We denote by $\gS(q,p)$ the point 
$(q, \tilde p)$, so that $\gS$ is a well defined involution of $T^*\breve M-\Gamma$. We denote by $\ms(q,p)>0$ the proportionality coefficient, defined by 
$$
\partial_pH(\gS(q,p))= -\ms(q,p)\partial_pH(q,p).
$$
We  study more precisely $\ms$ and $\gS$ in Section \ref{sec-sym}, where we prove:

\begin{prop}\label{prop-sym}
	The function $\ms$ extends with the values $1$ on $\Gamma$ to a locally Lipschitz function on $T^*\breve M$.
	
	The involution $\gS$ is smooth on $T^*\breve M-\Gamma$ and it extends to a $C^1$ fiber preserving diffeomorphism of $T^*\breve M$, which fixes $\Gamma$. We have, for $x=(q,\wp(q))\in \Gamma$,
	$$
	d\gS_{x}=\begin{bmatrix}Id &0\\2d\wp_q & -Id\end{bmatrix}.
	$$
\end{prop}

When $H$ is fixed and $x$ is a regular point of $H$, we call linearized energy level the kernel of the differential of $H$ at $x$, it is a hypersurface of $T_xT^*M$. It contains the vectorfield $V_H(x)$. The quotient of the linearized energy level by the vectorfield is the symplectic reduction of the linearized energy level, it has dimension $2d$ if $M$ has dimension $d+1$, and it
has a natural symplectic form descending from the one on $T^*M$. We call it the reduced tangent space.
Note that it does depend on $H$, and even on $u$. The reduced tangent spaces of a given Hamiltonian form a smooth vector bundle above the set of regular points of $H$. A hyperplane of $T_xT^*M$ transverse to the vectorfield at $x$ is called a linear section. It has a $2d$ dimensional intersection with the energy level, called a restricted linear section.
The restrictions of the symplectic form to restricted linear sections are symplectic forms.
The restricted linear sections at a point are symplectically isomorphic to the reduced tangent space.
Given two points on a same non-constant orbit, the differential of the flow defines a linear map between the tangent spaces at these points. This differential is symplectic and it preserves the linearized energy level, hence it gives rise to a symplectic linear map between the reduced tangent spaces, called the reduced linearized flow and denoted by $\Phi_s^t$.
If transverse sections are given near the points under consideration, then there is a well defined transition map between theses 
sections. The differential of this map preserves the linearized energy level, and its restriction is a symplectic  linear map between the restricted linear sections, called the restricted linearized transition map.
It  is conjugated to the reduced linearized flow.

The symmetry $\gS$ is preserving the energy $H$, hence its differential $d\gS_x$ maps the linearized energy level $\{dH_x=0\}$ at $x$ to the linearized energy level at $\gS(x)$. 
At points $x$ such that the  Hamiltonian vectors $V_H(\gS(x))$ is proportional to  $d\gS_x \cdot V_H(x)$,
(then $
V_H(\gS(x))=-\ms(x) d\gS_x\cdot V_H(x)$), 
the linear map $d\gS_x$ descends to a linear  map $\gR_x$ between the reduced energy spaces,
 called the reduced linearized symmetry. 
This holds at each points of round trip orbits, and at each points of $\Gamma$.
In the reversible case, the reduced linearized symmetry is antisymplectic at each point,
meaning that the pull-back of the natural symplectic form is minus the natural symplectic form.
 But in general, the reduced linearized symmetry is not necessarily antisymplectic at the points where it is defined.

\begin{defn}
	Let $\theta$ be a round trip orbit for $H$ of minimal period $T$ and les $\nu_0\in \Rm$ be a time such that $\theta(\nu_0) \in \Gamma$. We say that $\theta$ is reversible if :
	\begin{enumerate}
		\item
		The time $\nu_1:=\nu_0+T/2$ also satisfies $\theta(\nu_1)\in \Gamma$.
		\item The reduced linearized symmetries $\gR_{\theta(\nu_0)}$ and $\gR_{\theta(\nu_1)}$ are antisymplectic.
\item The reversibility identity 
		$$
		\gR_{\theta(\nu_0)}\circ \Phi^{\nu_0+T}_{\nu_1}\circ \gR_{\theta(\nu_1)}\circ\Phi_{\nu_0}^{\nu_1}=Id
		$$
		holds, where  $\Phi^t_s$ is the reduction of the linearized Hamiltonian flow  
		$\partial_x \varphi(t-s,\theta(s))$.
			\end{enumerate}
	\end{defn}

If $H$ is reversible, then all round trip orbits of $H$ are  reversible, so the following result covers the case of reversible Hamiltonians.

\begin{thm}\label{thm-g4}
	Given a smooth convex Hamiltonian $H: T^*M\rightarrow \Rm$, 
	the following property is satisfied for generic $u\in C^{\infty}(M)$ :

	All \textbf{reversible} periodic 
	orbits of $H+u$ which have zero energy  have their restricted linearized 
	return map in the complement of $\Upsilon$.
\end{thm} 

The proof consists in using the reversibility relation to prove a perturbation Lemma, see Section \ref{sec-rmr}.
We finally make  precise the intuition that round trip orbits should exist only in specific situations where some kind of reversibility holds :

\begin{thm}\label{thm-g5}
	Given a smooth convex Hamiltonian $H: T^*M\rightarrow \Rm$, 
	the following property is satisfied for generic $u\in C^{\infty}(M)$ :

	All \textbf{round trip} periodic
	orbits of $H+u$  are reversible.
\end{thm} 

Theorem \ref{thm-g3} obviously follows from Theorem \ref{thm-g4} and \ref{thm-g5}.
Theorem \ref{thm-g5} is proved in Section \ref{sec-rev}. The rough idea is that when a round trip  orbit is not reversible, it can be turned to a neat orbit by adding a small potential.
In the course of the proofs, we will use the notion of a chord :

\begin{defn}
	A chord (of energy zero) for $H$  is a point $(t,x)\in ]0,\infty[\times \Gamma$ such that 
	$H(x)=0$ and $\varphi(t,x)\in \Gamma$. 
	
	The chord $(t,x)$ is said transverse if $x$ is not a fixed point of the flow, and if 
	the flow map $\varphi$, in restriction to 
	$\Rm\times \Gamma_0$, is transverse to $\Gamma$  at the point $(t,x)$, where $\Gamma_0=\Gamma\cap \{H=0\}$.

	The chord $(t,x)$ is called minimal if $(s,x)$ is not a chord for $s\in ]0,t[$.
\end{defn}

Note, in the above definition, that $\Gamma_0$ is a submanifold of $\Gamma$ near  $x$ if this point is not a singular point of
the Hamiltonian.

We know that each round trip orbit contains exactly two points of $\Gamma$, each of which is the starting point of a chord.
In the reversible case, the converse is also true : The starting point of a chord belongs to a round trip orbit.
If a chord is transverse, it persists after a small perturbation of the Hamiltonian. 
In the reversible case, if $H$ admits a transverse chord, then $H+u$ also admits a chord for small $u$, and $H+u$ is still reversible, so $H+u$ still admits a round trip orbit.  
The following result on the transversality of chords follows from classical  tools detailed in Section \ref{sec-tc}:

\begin{prop}\label{prop-gen-chord}
	In the context of Theorem \ref{thm-g1}, for generic $u$, all minimal chords of energy zero  are transverse.
\end{prop}

Let us now make some general comments concerning the terminology used in the sequel.
We shall consider potentials $u\in C^{\infty}(M)$ both as functions on $M$ and as functions on $T^*M$,  which means that we
still denote $u$ the function $u\circ \pi$, where $\pi :T^*M\lto M$ is the canonical projection.
For each $u$, we thus consider $H+u$ as a Hamiltonian on $T^*M$, we denote by $V_H(x,u)$ the corresponding vectorfield, and by
 by $\varphi(t,x,u)$ the flow. We want to think of these maps as smooth, but since $C^{\infty}(M)$ is not a Banach space, some care is useful. The expression  $\partial_u \varphi \cdot v$ shall be considered as a directional 
 derivative in the direction $v$. More generally, we will consider finite dimensional subspaces $E\subset C^{\infty}(M)$, 
 and the derivative in the direction $E$, which can be defined as the derivative $\partial_v \tilde \varphi(t,x,u,0)$
 of the modified map $\tilde \varphi(t,x,u,v):=\varphi(t,x,u+v)$. It is a well defined linear map from $E$ to the appropriate tangent space of $T^*M$, and it depends continuously on $(t,x,u)$.

\section{The Symmetry}\label{sec-sym}

We study the symmetry $\gS$ and the function $\ms$ defined in the introduction, and prove Proposition \ref{prop-sym}.
We assume, without loss of generality in this section, that $\breve M=M$.
We denote $\gS(q,p)=(q,\gs(q,p))$. 
We first prove the smoothness of $\ms$ and $\gS$ outside of $\Gamma$. The equations defining $\ms$ and $\gs$ are
$$
\partial_pH (q, \gs(q,p))+\ms(q,p)\partial_pH(q,p)=0\quad, \quad H(q,\gs(q,p))-H(q,p)=0.
$$
We can apply the implicit function theorem to these equations provided the matrix
$$
D(q,p,\gs,\ms):=\begin{bmatrix}
\partial^2_{pp}H(q, \gs) & \partial_pH(q, p)\\ (\partial_pH(q,\gs))^t& 0 
\end{bmatrix},
$$
is invertible at the point $(q,p,\gs,\ms)$.
Since $(q,p)\not \in \Gamma$, we have $(q,\gs)\not \in \Gamma$, the vectors  
$
\partial_pH(q,p)$ 
and 
$\partial_pH(q,\gs)$
are both non zero, and they are negatively proportional.
To verify that $D(q,p,\gs, \ms)$ is invertible we consider an element
 $(\rho, r)$ of its kernel. Recall that $B:= \partial^2_{pp}H(q,\gs)$ is positive definite. 
We have $B\rho+r\partial_pH(q,p)=0$, hence $\rho=-rB^{-1}\partial_pH(q,p)$.
The second equation 
$\langle \rho, \partial_p H(q,\gs)\rangle =0$ (scalar product),
is thus equivalent to 
 $$
 r\langle B^{-1}\partial_pH(q,p), \partial_pH(q,\gs)\rangle =0.
 $$
The scalar product is different from zero because $B^{-1}$ is positive definite and because the vectors $\partial_pH(q,p)$ and $\partial_pH(q,\gs)$ are proportional and not null. We deduce that $r=0$, and then that $\rho=0$.

Before studying  $\gS$ and $\ms$ near $\Gamma$, we need a Lemma.

\begin{lem}\label{lem-sym}
	Let $U$ be an open set in some  $\Rm^k$, and let $f(y,x):U\times \Rm^d\lto \Rm$ be a smooth ($C^3$ is enough) function such that 
	$f_y:x\lmto f(y,x)$ is minimal at $x=0$ and  convex with positive definite Hessian at each point.
	Then for each $y\in U, x\neq 0$, there exists a unique point $\tilde x \neq x$ which is proportional to $x$ and such that 
	$h(y,\tilde x)=h(y,x)$.
	The involution $S:(y,x)\lmto (y,\tilde x)$ is smooth outside of $U\times \{0\}$ and extends to a $C^1$ diffeomorphism $S$ of $U\times\Rm^d$ fixing $U\times \{0\}$
	and satisfying $dS_{(y,0)}=
	\begin{bmatrix} Id & 0\\ 0 &-Id\end{bmatrix}$.
\end{lem}

\proof
We assume, without loss of generality, that $f(y,0)=0$ (adding a function of $y$ does not change the involution).
We search the map $S$ under the form $S(y,x)=(y,-s(y,x)x)$, $s(y,x)>0$. 
It is clear that the map $S$ is well-defined outside of $0$. The implicit function Theorem implies that it is smooth away from $0$. Indeed,
we can set $g(s,x,y):= f(y,-sx)$ and see  that $\partial_sg=-\partial_xf_x\cdot x<0$ if $x\neq 0$,  by convexity. Hence locally the unique positive  solution $s$ of the equation $g(s,x,y)=f(y,x)$ depends smoothly on $(y,x)$. 
It is moreover clear that the map $S$ extends by continuity on $U\times\Rm^d$ with $S(y,0)=0$.

We claim that $s$ extends to a locally Lipschitz function on $U\times \Rm^d$, with the value $s(y,0)=1$.
Assuming the claim for the moment, we return to the map $S(y,x)=(y,-s(y,x)x)$. The differential at a point $x\neq 0$ is 
$$dS(y,x)=
\begin{bmatrix}
I&0\\ -\partial_ys(y,x) \otimes x & -\partial _x s(y,x)\otimes x-s(y,x) I
\end{bmatrix}.
$$
Here we denote by $l\otimes x$ the linear map $v\lmto l(v)x$ when $l$ is a linear form. Since the differential of $s$ is bounded, we see that 
$dS$ extends by continuity on $U\times \{0\}$, with the value 
$$dS_{(y,0)}=\begin{bmatrix} Id & 0\\ 0 &-Id\end{bmatrix}.$$
This implies that $S$ is actually $C^1$.

We now prove the claim on $s$. We work in the neighborhood of a given point of $U$, called $0$.
We denote $A(y):= \partial^2_{xx}f(y,0)$. We can suppose by possibly reducing $U$ to a smaller neighborhood of $0$ that $a I\leq A(y) \leq I/a$ for some $a>0$ and for all $y$.

Since the third derivative of $f$ is bounded near $(0,0)$, we have, locally
$$
a|x|^2/4\leq f\leq |x|^2/a
$$
and, since that $f(y,x)=f(y,-s(y,x)x)$,
$$
as^2(y,x)|x|^2/4\leq f\leq s^2(y,x)|x|2/a
$$
hence 
$$
a^2/4\leq s^2(y,x)\leq 4/a^2
\quad,\quad
a/2\leq s(y,x)\leq 2/a.
$$
In the next computations, we denote by $O(|x|^k)$ a function of $(y,x)$ which is bounded by $C|x|^k$, with a locally uniform 
constant $C$.
We have
$$
f(y,x)=Ax^2/2+O(|x|^3)
$$
which implies, in view of the bound already obtained on $s$, that
$$
f(y,x)=f(y,-s(y,x)x)=s^2(y,x)Ax^2/2+O(|x|^3).
$$
These equalities imply that 
$$
(s^2(y,x)-1)Ax^2/2=O(|x|^3),
$$
hence that 
$
s^2(y,x)-1=O(|x|)
$
and finally,
$$
s(y,x)=1+O(|x|).
$$
As a consequence,  $S(y,x)=(y,-x)+O(|x|^2)$. 
Differentiating the equation 
$f(y,-s(y,x)x)=f(y,x)$
with respect to $x$ at $x\neq 0$ gives:
$$
-\big(\partial_xf(S(y,x))\cdot x \big) \partial_x s(y,x)-s(y,x) \partial_x f(S(y,x))=\partial_x f(y,x).
$$
Since $
\partial_x f(y,x)=A(y)x+O(|x|^2)
$
and $s=1+O(|x|)$, we obtain
$$(\partial_xf({S(y,x)})\cdot x ) \partial_x s(y,x)=O(|x|^2).
$$
We deduce that $\partial_x s$ is bounded near $(0,0)$ observing that 
$\partial_xf({S(y,x)})\cdot x=A(y)x^2+O(|x|^3)\geq a|x|^2/2$,

Similarly, differentiating the equation defining $s$ with respect to $y$ yields
$$
-\big(\partial_xf(S(y,x))\cdot x \big) \partial_y s(y,x)+ \partial_y f(S(y,x))=\partial_y f(y,x).
$$
Moreover, from $\partial_xf(y,0)\equiv 0$, we deduce that $\partial^2_{yx}f (y,0)=0$, hence that 
$\partial_yf(y,x)= O(|x|^2)$.
 We obtain that 
$$(\partial_xf({S(y,x)})\cdot x ) \partial_y s(y,x)=O(|x|^2)
$$
and conclude as above that $\partial_y s$ is locally bounded.
\qed

The map 
$$\mathfrak{h}: (q,p)\lmto (q, \partial_pH(q,p))
$$
is a  diffeomorphism from $T^*M$ to an open neighborhood of the zero section in  $TM$ (this is the legendre transform, it is onto if one adds the assumption that $H$ be superlinear, which is not useful here).
It maps $\Gamma$ to the zero section.  We denote its inverse by 
$
\mathfrak{g}
$ and denote $f:= H\circ \mathfrak{g}$. This is the energy expressed in the tangent bundle
(it is not the Lagrangian).
We claim that $\partial^2_{vv}f$ is positive definite on the zero section.
Indeed, let us fix the first coordinate $q$ for the moment, and denote by $h$ the map $p\lmto \partial_pH(q,p)$ and by $g$  its inverse, so that $\mathfrak{h}(q,p)=(q,h(p)),\mathfrak{g}(q,v)=(q,g(v))$.
We have
$$
\partial_vf_{(q,v)}=
h(g(v))\circ \partial_vg_v
=v\cdot \partial_vg_v
$$
and $\partial_vg_v= (\partial_ph_{g(v)})^{-1}=\big(\partial^2_{pp}H_{\mathfrak{g}(q,v)}\big)^{-1}$.
At $v=0$, we have
$$
\partial^2_{vv}f_{(q,0)}=\partial_vg_0=\big(\partial^2_{pp}H_{\mathfrak{g}(q,0)}\big)^{-1},
$$
and this matrix is positive definite.

We apply Lemma \ref{lem-sym}  to $f=H\circ \mathfrak{g}$ and deduce the existence of a unique positive function $s(q,v)$,
locally Lipschitz near the zero section, such that 
$$H\circ \mathfrak{g} (q, -s(q,v)v)=H\circ \mathfrak{g}(q, v).
$$
Moreover, the map $S(q,v):= (q, -s(q,v)v)$ is $C^1$ with
$dS(q,0)=
\begin{bmatrix}
I&0\\0&-I
\end{bmatrix}.$

Returning to $\ms$ and $\gS$, the equation 
$
\partial_pH(\gS(q,p))= -\ms(q,p)\partial_pH(q,p)
$
can be rewritten
$$
\gS(q,p)=\mathfrak{g} (q,-\ms(q,p) \partial_p H(q,p)),
$$
and since $H\circ \gS=H$ we obtain 
$$
H\circ \mathfrak{g} (q, -\ms(q,p) \partial_p H(q,p))
=H(q,p)=H\circ \mathfrak{g} (q, \partial_pH(q,p)).
$$
This equation implies that $\ms=s\circ \mathfrak{h}$, and then
$$
\mathfrak{h}\circ \gS(q,p)= (q, -\ms(q,p)\partial_pH(q,p)))
(q, -s(q,\partial_pH(q,p))\partial_pH(q,p)= S\circ \mathfrak{h}(q,p).
$$
Since
$$
\ms=s \circ \mathfrak{h}
\quad,\quad \gS=\mathfrak{g}\circ S\circ \mathfrak{h},
$$
 the regularity claimed for $\ms,\gS$ in the neighborhood of $\Gamma$ follows from the regularity of $s,S$ in the neighborhood of the zero section.
Let us finally compute the differential of $\gS$ at a point $x=(q, \wp(q))$ of $\Gamma$. For such a point, $\mathfrak{h}(x)$
belongs to the zero section, and
$S\circ \mathfrak{h}(x)=\mathfrak{h}(x)$, hence 
$$
d\gS_x
=d\mathfrak{g}_{\mathfrak{h}(x)}\circ dS_{\mathfrak{h}(x)}\circ d\mathfrak{h}_x
=(d\mathfrak{h}_x )^{-1}\circ dS_{\mathfrak{h}(x)}\circ d\mathfrak{h}_x.
$$
Taking coordinates, we have, at $x$, 
$$
d\mathfrak{h}=
\begin{bmatrix}
	I&0\\ \partial^2_{qp}H &\partial^2_{pp}H 
\end{bmatrix}
\quad, \quad
(d\mathfrak{h})^{-1}=
\begin{bmatrix}	
	I&0\\ -(\partial^2_{pp}H)^{-1} \partial^2_{qp}H &(\partial^2_{pp}H)^{-1} 
\end{bmatrix}
\quad, \quad
dS_{\mathfrak{h}(x)}=
	\begin{bmatrix}
I&0\\0&-I
	\end{bmatrix}
$$
from which follows that 
$$
d\gS_x=\begin{bmatrix}	
	I&0\\ -2(\partial^2_{pp}H_x)^{-1} \partial^2_{qp}H_x &-I
\end{bmatrix},
$$
and we can check by differentiating the equality 
$
\partial_pH(q, \wp (q))=0
$
 that 
$$
d\wp_q=-
(\partial^2_{pp}H_{(q, \wp(q))})^{-1}
\partial^2_{qp}H_{(q, \wp(q))}.
$$
That the lower right block is $2d\wp_q$ can be recovered also from the fact that $\gS$ is fixing $\Gamma$.
\qed

\section{Projected orbits}\label{sec-proj}

We study the projection of orbits and prove Proposition \ref{prop-neat}.
Let $\theta(t)=(Q(t),P(t))$ be a periodic orbit of minimal period $T$ and energy $e$.
We say that $s\in \Rm/T\Zm$ is a neat time if $\dot Q(s) \neq 0$ and if there exists no time $\sigma \neq s$ in $\Rm/T\Zm$ such that $Q(s)=Q(\sigma)$.
We say that $s$ is a degenerate time if $\dot Q(s)=0$, and we say that $s$ is a self-intersection time if there exists $\sigma \neq s$ in $\Rm/T\Zm$ such that $Q(s)=Q(t)$.

\begin{lem}
	There are finitely many degenerate times, and they are not self-intersection times. 
\end{lem}

\proof
If $s$ is a degenerate time, then $P(s)=\wp(Q(s))$ is the only minimum of the function $p\lmto H(Q(s),p)$ on $T^*_{Q(s)}M$, hence there is no other point $P\neq P(s)$ such that 
$H(Q(s),P)=H(Q(s), P(s))$. This implies that $s$ is not a self-intersection time.

Since $\dot Q(s)=0$ and since the orbit is not a fixed point,  $\dot P(s)\neq 0$.
Differentiating at $t=s$ the equality $\dot Q(t)=\partial_p H(Q(t), P(t))$ gives 
$$
\ddot Q(s)=\partial^2_{pp}H(Q(s), P(s))\dot P(s)\neq 0
$$
and we deduce that the degenerate time $s$ is isolated.
\qed

The next Lemma implies that an orbit which admits a neat time is neat.

\begin{lem}
	The set of neat times is open in $\Rm/T\Zm$.
	\end{lem}

\proof
Let $s_n\lto s, s_n\neq s$ be a sequence of times which are not neat. We will prove that the limit $s$ is not neat (it's either degenerate or self-intersection). Since there are finitely many degenerate times, we can asume by taking a subsequence that $s_n$ are self intersection times. 
Let $\sigma _n\neq s_n$ be such that $Q(s_n)=Q(\sigma_n)$.
Up to taking a subsequence, we can assume that $\sigma_n$ has a limit $\sigma$, and then $Q(\sigma)=Q(s)$. If $\sigma\neq s$, then $s$ is a self-intersection time.
If $\sigma= s$, then the equality $Q(s_n)=Q(\sigma_n)$ implies that $Q$ is not one to one near $t=s$, hence that $\dot Q(s)=0$.
\qed

We say that $s$ is a transverse time of self-intersection if for each  $\sigma\neq s \mod T$ such that $Q(\sigma)=Q(s)$ the vectors
$\dot Q(\sigma)$ and $\dot Q(s)$ are linearly independant.
If $s$ is a non-transverse self-intersection time, there exists exactly one $\sigma \neq s\mod T$ such that 
$\theta(\sigma)=\gS (\theta(s))$.

\begin{lem}\label{lemTSI}
	If a self-intersection time is transverse, then it is an isolated self-intersection time.
\end{lem}

\proof
Consider a transverse self-intersection time $s$, and a sequence $s_n\lto s, s_n\neq s$ of self-intersection times. For each $n$, there exists $\sigma _n \neq s_n$ such that 
$Q(\sigma _n)=Q(s_n)$. By taking a subsequence, we can assume that $\sigma _n$ has a limit $\sigma $, and then $Q(\sigma )=Q(s)$.
Since $s$ is not a degenerate time we must have $\sigma\neq s$, and the vectors  $\dot Q(s)$ and $\dot Q(\sigma )$ are linearly independent (because $s$ is assumed to be a transverse self-intersection time).
This implies that the geometric curves $\{Q(t), |t-\sigma|<\epsilon \}$ and $\{Q(t), |t-s|<\epsilon \}$ intersect only at the point $Q(s)=Q(\sigma)$ when $\epsilon>0$ si small. This is in contradiction with the existence of the sequences $s_n$ and $\sigma_n$ such that $s_n\lto s$, $\sigma_n\lto \sigma$ and $Q(s_n)=Q(\sigma_n)$. 
\qed

From now on, we assume that $\theta$ does not admit a neat time, and prove that it is a round trip orbit.
Let $R\subset \Rm/T\Zm$ be the set of non-degenerate times.
In the absence of neat points, all points of $R$ are self-intersection points. Since $R$ has a finite complement, it contains no isolated point hence by Lemma \ref{lemTSI}, there is no tranverse self-intersection; all points of $R$ are
 non-transverse self intersections.
 For each $s\in R$, there exists a unique time $\sigma(s)$ such that $\theta(\sigma(s))=\gS(\theta(s))$. 
If $\xi$ is smooth function on $M$ (a coordinate) such that $d\xi_{Q(s)}\cdot \dot Q(s)\neq 0$, then the equation $Q(\sigma(s))=Q(s)$ implies that $\xi\circ Q(\sigma(s))=\xi\circ Q(s)$.
By the implicit function theorem, the solution $\sigma(s)$ of this equation is smooth and locally decreasing. The function $\sigma$ is thus smooth 
and decreasing on each connected component of $R$.
 It is an involutive diffeomorphism of $R$,
which has no fixed point. Since all orientation-reversing homeomophism of the circle or of the interval  have a fixed point, we deduce that there are at least two points 
in the complement of $R$.
We consider a maximal interval $]\nu_0,\nu_1[$ in $R$, with boundaries $\nu_0\neq \nu_1$. Since $Q(\sigma(\nu_1))=Q(\nu_1)$, and since $\nu_1$ is a degenerate time (hence not a self-intersection time),
 we deduce that $\sigma(\nu_1)=\nu_1$, and similarly $\sigma(\nu_0)=\nu_0$. This implies that $\sigma$ extends as a homemorphism which maps
 $]\nu_0, \nu_1[$ onto $\Rm/T\Zm-[\nu_0, \nu_1]$.
 
 The last step is to prove that $\sigma$ is actually smooth around $\nu_i$, with differential equal to $-1$. We assume without loss of generality that $0$ is one of the fixed points of $\sigma$. Recall that this is a degenerate point of $Q$, so $Q'(0)=0$, $Q''(0)\neq 0$. We pick a smooth function $\xi$ on $M$ such that $\xi\circ Q(0)=0$ and  $(\xi\circ Q)''(0)>0$. We denote by $f:\Rm\lto \Rm $ the function $\xi\circ Q$.
 Recall that a smooth function $f:\Rm\lto \Rm$ satisfying $f(0)=0$ can be written $f(t)=tf_1(t)$ with a smooth function $f_1$ (
$f_1(t)=\int_0^1 f'(ts)ds$, it  has one derivative less than $f$ at $t=0$). Applying this twice yields that $f(t)=t^2 f_2(t)$, with $f_2(0)> 0$. Taking the square root, we obtain that $f(t)=g(t)^2$, with $g=t\sqrt{f_2}$ smooth and $g'(0)>0$.
 The equation $Q(\sigma(t))=Q(t)$ is now equivalent to $g (\sigma(t))=-g(t)$. The implicit function theorem can be applied to this equation and yields that $\sigma$ is smooth at $t=0$ (it actually has two derivatives less than $Q$) with $\sigma'(0)=-1$.
  \qed

Note that 
$
\sigma'(t)=-\ms(\theta(t)),
$
which also implies that $\sigma$ is $C^1$, with $\sigma'(\nu_i)=-1$. The above proof provides a better regularity for $\sigma$.
We finish with a remark on multiple intersections. We say that $t$ is a multiple intersection time if there are at least three different times in $\Rm/T\Zm$  (including $t$) at which $Q$ takes the value $Q(t)$.

\begin{lem}\label{lem-multi}
Each periodic orbit $\theta$ has finitely many multiple intersections.
\end{lem}

\proof
If there was infinitely many multiple intersections, there would exist injective sequences $s_n\lto s$, $t_n\lto t$, $\tau_n \lto \tau$
such that $Q(s_n)=Q(t_n)=Q(\tau_n)$, and such that the three times $s_n,t_n,\tau_n$ are distinct for all $n$.
At the limit, $Q(s)=Q(t)=Q(\tau)$. 

If two of the limits $s,t,\tau$ are equal, say $s$ and $t$, then, $s$ is a degenerate point (because the curve $Q$ is one to one
near a non-degenerate point). Since a degenerate point  can't be an intersection point, we deduce that  $\tau=s=t$.
But since $Q''(s)\neq 0$, the curve $Q$ takes each value at most twice near $s$, a contradiction.

The second possibility is that $s, t, \tau$ are distinct, and non-degenerate. But then at least two of the three derivatives 
$Q'(s), Q'(t), Q'(\tau)$ are linearly independant, say $Q'(s)$ and $Q'(t)$. This is in contradiction with the existence of the intersection points $s_n\lto s$ and $t_n\lto t$.
\qed

\section{Reversible symplectic matrices}

We consider here the group $\Sp(2d)$ of real symplectic matrices and the space $\A(2d)$ of antisymplectic involutions, which are $2d\times 2d$
 matrices $R$ such that $R^2=Id$ and $R^*\omega=-\omega$, where $\omega$ is the standard symplectic form on $\Rm^{2d}$. In matrix form,
 this second equation can be rewritten $R^tJR=-J$, with the usual symplectic matrix
 $$J=
 \begin{bmatrix} 0&I\\-I&0
 \end{bmatrix}.
 $$
 Two examples are  the map $R_0:(q,p)\lmto (q,-p)$ and $R_1:(q,p)\lmto (p,q)$. We also denote by $\M(2d)$ the space of square real matrices of size $2d$. The content of this section is partly inspired by \cite[Section V]{D}, but Devaney studies there the space of all  reversible matrices, not the anti-symplectic reversible ones.

 \begin{prop}\label{prop-conj}
 	The subset $\A(2d)\subset \M(2d)$ is a connected algebraic submanifold without singularity of dimension $d(d+1)$.
 	The action of $\Sp(2d)$ on $\A(2d)$ by conjugacy is transitive and submersive, meaning that the map
 	$$
 	\Sp(2d)\ni M\lmto M^{-1}RM\in \A(2d)
 	$$
 	is a surjective submersion for each $R\in \A(2d)$.
 \end{prop}
 
 \proof 
 The eigenspaces $E_1(R)$ and $E_{-1}(R)$ associated to 
 an element $R$ of $\A(2d)$ are isotropic, and since $\Rm^{2d}=E_1(R)\oplus E_{-1}(R)$,
 they must each have dimension $d$ and be Lagrangian.
 Conversely,  two transverse Lagrangian subspaces  of $\Rm^{2d}$ determine a unique element of $\A(2d)$. Since the group  $\Sp(2d)$ acts transitively on ordered pairs of transverse Lagrangian subspaces
 (see e.g. \cite[Theorem 1.26]{dG}), it also acts transitively by conjugacy on $\A(2d)$. The isotropy subgroup of  an element $R\subset \A(2d)$ is the set of symplectic isomorphisms which preserve each of the two Lagrangian spaces  $E_1(R)$ and $E_{-1}(R)$, it has dimension $d^2$. As a consequence, the rank of the map $\Sp(2d)\ni M\lmto M^{-1}RM\in \M(2d)$ is $\dim \Sp(2d)-d^2=d(d+1)$. 
 
 On the other hand, we can linearize the equations $R^2=I$ and $R^TJR=-J$ which determine $\A(2d)$ at the point 
 $R_0=\begin{bmatrix} I&0\\0&-I\end{bmatrix}$.
 We obtain the equations $R_0R+RR_0=0$ and $R_0JR+R^TJR_0=0$.
  A simple computation in block form reveals that the kernel of these linearized equations is the space of matrices of the form 
 $\begin{bmatrix} 0&A\\B&0\end{bmatrix}$ with $A$ and $B$ symmetric.
 This is a $d(d+1)$ dimensional linear subspace of $\M(2d)$. This implies that the orbit $\A(2d)$ of the action of $\Sp(2d)$ is an embedded submanifold of dimension $d(d+1)$.
 \qed

Given  $R\in \A(2d)$, we say that $M\in \Sp(2d)$ is $R$-reversible if $RM$ is an involution, and then $RM\in \A(2d)$. The set of $R$-reversible matrices is thus $R\A(2d)$, it is a submanifold of $\M(2d)$. The transitive action of $\Sp(2d)$ on $\A(2d)$
gives rise to a transitive action on $R\A(2d)$, given by 
$$
\Sp(2d)\times R\A(2d) \ni(M, L)\lmto RM^{-1}RLM\in R\A(2d).
$$
For each fixed $R$, this map is a submersion from $\Sp(2d)$ onto $R\A(2d)$.
If $R$ and $\tilde R$ are two elements of $\A(2d)$, then the submanifolds $R\A(2d)$ and $\tilde R\A(2d)$ are conjugated inside $\Sp(2d)$.

As is well known, the symplectic matrices which preserve the first component are of the form $\begin{bmatrix} I&0\\S&I\end{bmatrix}$ for some symmetric matrix $S$. Similaly, 
the elements of $\A(2d)$
which preserve the first component  are of the form $\begin{bmatrix} I&0\\S&-I\end{bmatrix}$ for some symmetric matrix $S$.

\begin{prop}\label{prop-gen}
	For each $R\in \A(2d)$, the matrices with multiple eigenvalues form a closed and nowhere dense set in $R\A(2d)$. So do the matrices having a fixed complex number $\lambda$ as an eigenvalue.
\end{prop}

\proof
We consider the algebric map $\Delta : R\A\lmto \Rm$ which, to each $M\in R\A(2d)$ associates the discriminant of its caracteristic polynomial. The matrix $X$ has a multiple eigenvalue if and only if $\Delta(X)=0$.
By  analytic continuity, either $\Delta$ is identically vanishing on $R\A(2d)$ or the set $\Delta^{-1}(0)$ has empty interior.
So we just need to show the existence of an element of $R\A(2d)$ without multiple eigenvalues.
We consider $M\in \Sp(2d)$ such that 
$M^{-1}RM=R_1=
 \begin{bmatrix} 0&I\\I&0\end{bmatrix}
 $.
Then, we consider a diagonal matrix $X$, with diagonal elements $(x_1, \ldots, x_d, 1/x_1, \ldots, 1/x_d)$,
with $1<x_1< \cdots <x_d$. This matrix $X$ is symplectic and $R_1$-reversible, and it has simple eigenvalues.
Then, $MXM^{-1}$ is symplectic, it has simple eigenvalues, and it is $R$-reversible since
$$
MXM^{-1}RMXM^{-1}R=MXR_1XR_1M^{-1}=MM^{-1}=I.
$$
We conclude similarly concerning the matrices having $\lambda$ as an eigenvalue by considering the algebraic map $M\lmto \det (M-\lambda I)$. 
\qed
To sum up:

\begin{cor} \label{cor-RA}
	Given  $R\in \A(2d)$, the intersection $\Upsilon\cap R\A(2d)$ is an $F_{\sigma}$ with  empty interior in $R\A(2d)$.
	\end{cor}

\section{Transverse Chords}\label{sec-tc}
We prove Proposition \ref{prop-gen-chord} using a variation on  the parametric transversality principle, see \cite{A63}.
We  represent chords by  the set  $\mC\subset ]0,\infty[\times \Gamma\times C^{\infty}(M)$  of triples $(t,x,u)$ 
such that $(t,x)$ is a  chord (of energy $0$) for $H+u$. It is a closed subset of the product. The subset $\mC^n\subset\mC$ of chords which are not transversal is also closed.

\begin{lem}\label{lem-open}
	The set  of minimal chords is open in $\mC$.
\end{lem}

\proof
Let $(T_n,x_n,u_n)\lto (T,x,u)$ be a converging sequence of non-minimal chords. There exist times
$S_n\in ]0, T_n  [$ such that $(S_n,x_n,u_n)$ are also chords, and, up to a subsequence, 
$S_n$ has a limit  $S$ in $[0,T]$. 
If $x$ is a fixed point, then $(T,x,u)$ is not minimal.
If $x$ is not a fixed point, 
the vector-field 
$V_{H+u}$ is vertical at $x$, hence not tangent to $\Gamma$. Then there exists $\epsilon>0$ such that   $\varphi(t, x_n,u_n)\not \in \Gamma$ 
for  $t\in ]0, \epsilon[$ provided $n$ is large enough. This implies that $S>0$.
Similarly, $V_{H+u}$ is vertical at $\varphi(T,x,u)$ and this implies that $S<T$.
Finally, we have $S\in ]0,T[$, and $(S,x,u)$ is a chord, hence $(T,x,u)$ is not minimal.
\qed

The following local version of Proposition \ref{prop-gen-chord} implies Proposition \ref{prop-gen-chord}:

\begin{prop}\label{prop-loc}
	Let $(T,\theta)$ be a minimal chord of $H$. Then there exists an open neighborhood $\mC_{loc}$ of $(T,\theta,0)$ in 
	$\mC$ such that the projection of $\mC^n_{loc}:=\mC^n\cap\mC_{loc}$ on $C^{\infty}(M)$ has empty interior.
\end{prop}

Let us first explain how Proposition \ref{prop-loc} implies Proposition \ref{prop-gen-chord}.
Let $\mC'$ be the set of chords which are minimal and not transverse. 
 For each $(t,x,u)\in \mC'$, we can apply Proposition \ref{prop-loc} to $H+u$ at point $(t,x)$ 
and get the existence of an open neighborhood
$\mC'_{loc}$ of $(t, x, u)$ in $\mC'$ whose projection on $C^{\infty}(M)$ has empty interior.
Moreover, $\mC'_{loc}$ is locally closed (the intersection of an open set and of a closed set), hence it is an $F_{\sigma}$ in 
the product $]0,\infty[\times \Gamma\times C^{\infty}(M)$.
Recalling that the projection of a closed set of the product on the last factor is an $F_{\sigma}$, we deduce that the projection of $\mC'_{loc}$ is an $F_{\sigma}$ with empty interior.
The separable metric space $\mC'$ can be covered by countably many open neighborhoods $\mC'_{loc}$ having this property, so its projection is a countable union of $F_{\sigma}$ with empty interior. By the Baire property, the projection of $\mC'$ is an
$F_{\sigma}$ with empty interior,   Proposition \ref{prop-gen-chord} is proved. 
\qed

\textsc{Proof of Proposition \ref{prop-loc}.}
We fix a minimal chord $(T,\theta )$ and denote by $\theta(t)=(Q(t), P(t))$ the $H$-orbit of the point $\theta$.
Since $T>0$ and the chord is minimal, the orbit $\theta(t)$ is not constant, hence, locally near $\theta(0)$ and near $\theta(T)$, the energy level is a submanifold, 
and it is 
transverse to $\Gamma$. 
The following key observation seems to first appear in \cite{A82}, Lemma 2:

\begin{lem}
	There exists a finite dimensional subspace $E\subset C^{\infty}(M)$, formed by potentials vanishing near $Q(0)$ and near $Q(T)$, 
	such that the Gâteau differential  
	$\partial_{t,u}\varphi(T,\theta, 0)
	$
	maps $\Rm\times E$ onto the linearized energy level of $H$ at 
	$\theta(T)$.
\end{lem}

We omit the proof, which is similar to \cite[Lemma 2]{A82} or \cite[Lemma 7]{AB2}.
Since the energy level of $H$ at $\theta(T)$ is transverse to $\Gamma$, this implies that $\varphi$, in restriction 
to $]0, \infty[\times \{\theta\}\times E$,  is transverse to $\Gamma$ at $(T,\theta,0)$.

Then, the  map 
$$
(t,v)\lmto \varphi(t, \theta, u+v) 
$$
is transverse to $\Gamma$ on $]0, \infty[_{loc}\times E_{loc}$ 
provided $(\theta,u)\in \Gamma_{loc}\times C^{\infty}_{loc}$, where $]0, \infty[_{loc}, \Gamma_{loc},C^{\infty}_{loc}, E_{loc}$ 
are sufficiently small open neighborhoods of $T,\theta, 0,0$ in the corresponding spaces.
We choose the neighborhoods such that the elements of $E$ vanish on $\Gamma_{loc}$ and such that  the restriction of
$ H+u$  to $\Gamma_{loc}$ has no critical point  for $u\in C^{\infty}_{loc}$. 

From now on, we fix $u\in C^{\infty}_{loc}$ and prove that $u$ does not belong to the interior of the projection 
of $\mC^n_{loc}:=\mC^n \cap \mC_{loc}$, with  $\mC_{loc}:=\mC\cap \big(]0, \infty[_{loc}\times \Gamma_{loc}\times C^{\infty}_{loc}\big)$.
We denote by $\Gamma^u_{loc}$ the submanifold
of $\Gamma_{loc}$ of equation $H+u=0$, or equivalently of equation $H+u+v=0$ for each $v\in E$.
The map
$$
\psi: ]0, \infty[_{loc}\times \Gamma^u_{loc}\times E_{loc}\ni (t,x,v)\lmto \varphi(t, x, u+v) \in T^*M
$$
is transverse to $\Gamma$. 
As a consequence, the set
$\psi^{-1}(\Gamma)$ is a smooth submanifold. 
The point $(t,x,v)$ belongs to $\psi^{-1}(\Gamma)$ if and only if $(t,x,u+v)$ is a chord of energy $0$.
Moreover, this chord is transverse if  $(t,x,v)$ is a regular point of the restriction $\pi_{rest}$ to $\psi^{-1}(\Gamma)$ 
of the projection to the third factor.
In order to check this, we   denote by $F$ the product $\Rm\times T_x\Gamma^u_{loc}$,  by $V\subset T_{\psi(t,x,v)}T^*M$ the linear vertical,
and by $\Psi:F\times E\lto V$ the composition of the differential $d\psi_{(t,x,v)}$ and of the projection on $V$ parallel to $T_{\psi(t,x,v)}\Gamma$. This linear map is onto.
By definition, the chord under consideration is transverse if  $\Psi$ maps $F\times \{0\}$ onto $V$. It is an easy exercise of linear algebra 
to check that this is equivalent to the fact that the kernel  of $\Psi$, which is the tangent space to $\psi^{-1}(\Gamma)$, projects onto $E$. 

Finally, we apply Sard's Theorem to the restricted projection $\pi_{rest}$ and deduce   the existence of arbitrarily small regular values $v$ of this 
restricted projection. For such a $v$, all chords of energy zero  of $H+u+v$ in $]0, \infty[_{loc}\times M_{loc}$ are transverse, 
meaning that $u+v$ does not belong to the projection of $\mC^n_{loc}$. As a consequence, the function $u$ does not belong to the interior of this projection, and this holds for each $u\in C^{\infty}_{loc}$.
\qed

\section{The return map of reversible orbits}\label{sec-rmr}

We prove Theorem \ref{thm-g4}.
As in \cite{O08,RR}, we have to understand to what extent
the return map of a given orbit can be modified by adding a small potential to the Hamiltonian.
In that respect, the following result was established in \cite{RR} (completed in \cite{AB}):

\begin{prop}\label{prop-wo}
	Let $\theta(t)=(Q(t),P(t)), t\in [a,b]$ be an orbit segment of the convex Hamiltonian $H$. Assume that the first component
	$Q$ is an embedding of the interval $[a,b]$ into $M$, and let $U\subset M$ be an open set intersecting $Q(]a,b[)$.
	 Let $\mE\subset C^{\infty}(M)$ be the space of adapted potentials supported in $U$, more precisely the space of smooth potentials $u$ supported in $U$, null near $\{Q(a),Q(b)\}$, and satisfying $u=0$, $du=0$ on the image of $Q$ (which implies that $\theta$ is still an orbit segment for the potentials $H+u, u\in \mE$).
	 The map 
	 $$\mE \ni u\lmto \Phi_a^{b}(u)\in \Sp(2d)
	 $$
	  which, to each potential $u$ associates the restricted  linearized transition map for the Hamiltonian $H+u$ between times $a$ and $b$ is well defined once symplectic coordinates have been fixed, independently from $u$, on the reduced tangent spaces at $\theta(a)$ and $\theta(b)$ (these spaces do not depend on $u\in \mE$).
	  This map is weakly open, meaning that the image of each non-empty open set of $\mE$ contains a non-empty open set of $\Sp(2d)$.
	 \end{prop}

In this statement, we could replace "weakly open" by "open" using  \cite{LRR}. However, this is much harder to prove and not necessary for our present study.

When $\theta$ is a neat periodic orbit, one can find a time interval $[a,b]$ such that the above result applies, and such that moreover $Q([b,a+T])$ is disjoint from $Q(]a,b[)$, where $T$ is the minimal period of the orbit. Then, the above result can be applied with $U=M-Q([b,a+T])$, and the restricted return map $\Phi_a^{a+T}(u)$ of the periodic orbit $\theta$ for 
the Hamiltonian $H+u$  at $\theta(a)$  can be written
 $$\Phi_a^{a+T}(u)= \Phi_b^{a+T}\circ\Phi_a^{b}(u),
 $$
  where $\Phi_b^{a+T}$ does not depend on $u\in \mE$ because the elements $u\in \mE$ vanish near $Q([b,a+T])$. We deduce immediately that the map 
  $$\mE\ni u\lmto \Phi_a^{a+T}(u)$$
  is weakly open, which is a way of saying that we have sufficient possibility to change the return map by changing the potential.
  This is the main step in the proof of Theorem \ref{thm-g2}.

However, if the orbit $\theta$ is not neat,  it is not possible to chose the interval $[a,b]$ such that 
$Q([b,a+T])$ is disjoint from $Q(]a,b[)$, and then in the decomposition 
$\Phi_a^{a+T}(u)= \Phi_b^{a+T}(u)\circ\Phi_a^{b}(u),
$
both factors depend on $u$, and  it is not clear in general  how the composition behaves.
This is why the above discussion, which summarizes the strategy  of \cite{O08,RR} is not sufficient to obtain Theorem \ref{thm-g1}.

The proof of Theorem \ref{thm-g4} consists in solving this difficulty using the reversibility of the  orbit.
So we assume that  $\theta$ is reversible, and  also for definiteness that $\theta(0)\in \Gamma$
(then $\theta(T/2)\in \Gamma$).
We denote  by $\gR_t$ the reduction of the linearized symmetry  $d\gS_{\theta(t)}$.
In particular $\gR_0$ and $\gR_{T/2}$ are antisymplectic involutions  of the reduced tangent spaces at $\theta(0)$
and $\theta(T/2)$. These reduced tangent spaces can be represented by $2d$-dimensional subspaces transverse to the vectorfield
and contained in the linearized energy level. These transverse subspaces can be chosen invariant under the linearized symmetries
$d\gS_{\theta(0)}, d\gS_{\theta(T/2)} $, and then $\gR_0, \gR_{T/2}$ are just the restrictions of these linearized
symmetries to the linear sections.

   \begin{prop}\label{prop-cr}
   	Let $\theta$ be a reversible orbit of minimal period T, satisfying $\theta(0)\in \Gamma$. 
   	 Let  
   	 $\mE\subset C^{\infty}(M)$ be the space of potentials which satisfy $u=0, du=0$ on  $Q([0,T/2])$, and such that 
   	moreover $u\equiv 0$ in a neighborhood of $\{Q(0), Q(T/2)\}$ (we call them potentials adapted to the chord $(\theta(0), T/2)$).
   	Assume that there exists an open subset $\mE_{loc}\subset \mE$ such that the orbit $\theta$ is reversible for $H+u$ for each 
   	$u\in \mE_{loc}$.
   	Then the reduced return map $\Phi_0^T(u)$ along the orbit $\theta$ for $H+u$  belongs the space $\gR_0\A(2d)$ of $\gR_0$-reversible symplectic matrice for each 
   	$u\in  \mE_{loc}$, and the image $\Phi_0^T(\mE_{loc})$ contains an open subset of $\gR_0\A(2d)$.
   	\end{prop}
   	
   	\proof
   	We denote by $\Phi_s^t(u)$ the reduction of the linearized flow $\partial_x \varphi^{t-s}(\theta(s),u)$.
On  $\mE_{loc}$,  the reversibility relation
   	$$\gR_0\circ \Phi_{T/2}^T(u)\circ \gR_{T/2}\circ \Phi_0^{T/2}(u)=Id,
   	$$
 holds, it can be rewritten
   	$$
   	\Phi^{T}_{T/2}(u)=\gR_0\circ (\Phi_0^{T/2}(u))^{-1}\circ\gR_{T/2}.
   	$$
 Then, for $u\in \mE_{loc}$, 
   $$\Phi_0^T(u)= \Phi^{T}_{T/2}(u)\circ \Phi_0^{T/2}(u)
   =\gR_0\circ (\Phi_0^{T/2}(u))^{-1}\circ\gR_{T/2}\circ \Phi_0^{T/2}(u).
   $$
   So the set
   $ \Phi_0^T(\mE_{loc})
   $
   is the image
   of the set $\Phi_0^{T/2}(\mE_{loc})$, which contains an open subset of $\Sp(2d)$ by Proposition \ref{prop-wo},
    by the map
   $$\Sp(2d)\ni M\lmto \gR_0 M^{-1}\gR_{T/2} M \in \gR_0\A(2d).
   $$
   This map  is  open by Proposition \ref{prop-conj}, the statement follows.
   \qed

In view of  Proposition \ref{prop-gen-chord}, Theorem \ref{thm-g4} follows from :

\begin{prop}
	The following property is satisfied for generic $u\in C^{\infty}(M)$:
	
	For each $x\in \Gamma\cap\{H+u=0\}$ such that the $(H+u)$-orbit of $x$ is reversible periodic of minimal period $2T$, and such that
	the chord $(T,x)$ is transverse, the reduced return map belongs to the complement of $\Upsilon$.
\end{prop}
 
 \proof
Let $\mM\subset ]0,\infty[\times \Gamma\times C^{\infty}(M)$ be the set of
 minimal transverse chords of energy $0$.
 In view of the implicit function theorem (see Appendix \ref{ap}) and of
 Lemma \ref{lem-open}, the projection  from $\mM$ to $C^{\infty}(M)$ is a local homeomorphism. Let $\mP\subset \mM$ be the subset of transverse minimal chords  $(T,\theta, u)$
 associated to reversible periodic orbits with return map in $\Upsilon$. In other words, the chord 
  $(T,\theta,u)\in \mM$ belongs to $\mP$ if an only
  if :
 \begin{itemize}
 	\item $\varphi(2T,\theta, u)=\theta$,
 	\item The reduced symmetries $\gR_{(\theta,u)}$ and $\gR_{(\varphi(T,\theta,u),u)}$ are antisymplectic.
\item $\gR_{(\theta,u)}\circ \Phi^{2T}_T(\theta, u)\circ 
 	\gR_{(\varphi(T,\theta,u),u)}\circ \Phi_0^T(\theta, u)=Id$,
 	\item $\Phi _0^{2T}(\theta,u)\in \Upsilon$.
 \end{itemize}
In these equalities, we denote by $\gR_{(x,u)}$ the reduced linearized symmetry at point $x$ for the Hamiltonian $H+u$, and by 
$\Phi_s^t(x,u)$ the reduced linearized flow of $H+u$ along the orbit of $x$ from time $s$ to time $t$.
 We have  to prove that the projection of $\mP$ on the third factor is an $F_{\sigma}$
 with empty interior. 
Since the projection from $\mM$ to the third factor is a local homeomorphism, and since $\mM$ is a separable metric space, 
it is enough to prove that $\mP$ is 
 an $F_{\sigma}$ with empty interior in $\mM$. The first claim is clear, using that $\Upsilon$ is an $F_{\sigma}$.
To prove the second claim, we argue by contradiction and assume that the interior of $\mP$ contains a point $(T, \theta, u)$.
 
 As in Proposition \ref{prop-cr}, we denote by $\mE\subset C^{\infty}(M)$ the space of potentials $v$ which are adapted to the chord 
 $\theta$, meaning that $v=0, dv=0$ on $Q([0,T])$ (the first coordinate of the $(H+u)$-orbit of $\theta$), and that $v$ is vanishing near 
 $Q(0)$ and near $Q(T)$. Since $(T,\theta, u)$ is in the interior of $\mP$, there exists an open neighborhood $\mE_{loc}$ of $0$
 in $\mE$ such that $(t,\theta, u+v)\in \mP$ for $v\in \mE_{loc}$. We deduce
 from  Proposition \ref{prop-cr} 
 that the set
 $\Phi_0^{2T} (\mE_{loc})$ contains an open subset of  $\gR_{0}\A(2d)$.
 In view of Corollary \ref{cor-RA}
 this contradicts the fact that $\Phi_0^T(\mE_{loc})$ is contained in $\Upsilon$.
\qed

\section{Non homogeneous linear systems.}\label{sec-nhls}

In the present section, we expose a result on non-homogeneous linear systems which is the key step for the proof
of Theorem \ref{thm-g5}.
We consider two smooth curves of Hamiltonian $2d\times 2d$ matrices $L_t$ and $\tilde L_t$ defined on $[0,T]$
which we think as the linearized systems along the two branches of a  round trip orbit near a given non-degenerate point.
We denote 
$$
L_t=
\begin{bmatrix}
C_t^T & B_t\\-A_t &-C _t
\end{bmatrix}
\quad, \quad
\tilde L_t=
\begin{bmatrix}
\tilde C_t^T & \tilde B_t\\-\tilde A_t &-\tilde C_t
\end{bmatrix},
$$
 where $A_t, B_t, \tilde A_t, \tilde B_t$ are symmetric $d\times d$ matrices and $C_t, \tilde C_t$ are arbitrary 
matrices of the same size. We assume that $\tilde B_t$ is invertible for each $t$.
We also consider two non-vanishing smooth functions $a(t)$, $\tilde a(t)$ on $[0,T]$ and consider the inhomogeneous linear systems
\begin{align*}
x'(t) & = a(t)L_t x(t)+a(t)b(t)\\
\tilde x'(t)&=\tilde a(t)\tilde L_t \tilde x(t)+\tilde a(t) b(t)
\end{align*}
as well as their homogeneous counterparts. Note, and this is a key point, that the same curve $b$ appears in both systems.
We think of $b$ as a perturbation created by adding a potential, so it will take values in $\{0\}\times (\Rm^d)^*$.
We denote by $\Psi_s^t, \tilde \Psi_s^t$ the resolvants of the homogeneous systems $x'=aLx, \tilde x'= \tilde a \tilde L \tilde x$.
We recall the standard :

\begin{lem}\label{lem-cr}
	Given a curve $R_s$ of matrices, the conjugacy  relation
	$$
	R_r\circ \Psi_s^r= \tilde \Psi_s^r \circ R_s
	$$
	holds for each $r$  near $s$ if and only if the differential system
	$$
	R'_s +a(s) R_s L(s) =\tilde a(s)\tilde L(s) R_s
	$$
	holds for each $r$ near $s$.
\end{lem}


\proof
Setting $\Delta^r=R_{r}\circ \Psi_s^r -\tilde \Psi_{s}^{r}\circ R_s$, we have
\begin{align*}
	\partial_r \Delta^r
	&=(R_{r}' +a(r)R_{r}L(r))\Psi_s^r-\tilde a(r)\tilde L(r)\tilde \Psi_{s}^{r}(u)R_s\\
	&=\big(R_{r}' +a(r)R_{r}L(r)-\tilde a(r)\tilde L(r)R_{r}\big)\Psi_s^r(u)+\tilde a(r)\tilde L(r)\Delta^r.
\end{align*}
We see from this differential equation that $\Delta^r\equiv 0$ near $r=s$ if and only if the equation 
$$
R_{r}' +a(r)R_{r}L(r)-\tilde a(r)\tilde L(r)R_{r}=0
$$
holds for $r$ near $s$.
\qed

We are interested in particular in conjugacies preserving the first coordinates, \textit{i.e.} of the form 
$R_t = \begin{bmatrix} I&0\\ * &* \end{bmatrix}$. The first block line of the differential equation
$
R_{r}' +a(r)R_{r}L(r)-\tilde a(r)\tilde L(r)R_{r}=0
$
indicates that the only possible conjugacy of that form is
$$
\gR_{t}:=
\begin{bmatrix}
I&0\\
(\tilde a(t)\tilde B(t))^{-1}\big( a(t)C^T(t)-\tilde a(t) \tilde C^T(t)\big) &  (a(t)/\tilde a(t)) \tilde B^{-1}(t) B(t)
\end{bmatrix}.
$$

\begin{prop}\label{prop-3co}
	In the context described above, assume that $\tilde B(t)$ is invertible for each $t\in [0,T]$.
	Assume that, for each smooth curve $b(t):\Rm \lto \{0\}\times (\Rm^d)^*$, compactly supported in $[0,T]$, the solutions 
	$x$ and $\tilde x$ of the two non-homogeneous linear systems with the initial condition $x(0)=0=\tilde x(0)$ have the same projection on $\Rm^d\times \{0\}$.
	Then the following conditions are satisfied for each $t\in [0,T]$ :
	\begin{enumerate}
		\item $\tilde a/a$ is constant,	
		\item The matrices $\gR_t$ are conformally symplectic of factor $\tilde a/a$, meaning here that 
		$$\gR_t=\begin{bmatrix}
		I&0\\ \gs_t &(\tilde a /a) I
		\end{bmatrix}$$
		with $\gs_t$ symmetric.
		\item  The homogeneous systems $x'(t)=a(t)L_tx(t)$ and $\tilde x'(t)=\tilde a(t)\tilde L_t\tilde x(t)$ are conjugated by $\gR_t$.
	\end{enumerate}
	Conversely, if  conditions 2 and 3 are  satisfied, then  $\gR_t$ is conjugating the non-homogeneous equations $x'= aLx+ab$ and $\tilde x'=\tilde a \tilde L x+\tilde a b$ for each curve $b$ compactly supported in $]0,T[$ and taking  
	values in $\{0\}\times (\Rm^d)^*$, and in particular the solutions emanating from $0$ have the same projection on $\Rm^d\times \{0\}$.
\end{prop}

It might sound surprising to obtain a converse using only two of the three conclusions. It means that conclusions 2 and  3 actually imply conclusion 1. We will give a direct proof of this fact in the course of the following proof.

\noindent

\proof
We first prove the converse. 
If $x(t)$ solves the equation $x'=aLx+ab$, then $y(t):= \gR_tx(t)$ satisfies
$$
y'(t)=\gR_t'\gR_t^{-1} y(t)+ a(t)\gR_t L_t\gR_t^{-1}y(t) +a(t) \gR_t b(t).
$$
From the third assumed conclusion, we have 
$$
\gR_t'\gR_t^{-1} +a(t)\gR_t L_t\gR_t^{-1}=\tilde a (t)\tilde L_t.
$$
From the form of $\gR_t$ and $b$, we see that $\gR_t b(t)=(\tilde a(t)/a(t))b(t)$.
So $y$ solves  the system  $y'=\tilde a\tilde Ly+\tilde ab$, which was our claim.

\qed

Let us now explain why conditions 2 and 3 imply condition 1.
More generally, if the matrices $\gR_t$ are of the form 
$$
\gR_t=\begin{bmatrix}
I&0\\ \gs_t & \alpha I
\end{bmatrix}
\quad, \quad 
(\gR_t)^{-1}=\begin{bmatrix}
I&0\\ -\alpha^{-1}\gs_t & \alpha^{-1} I
\end{bmatrix}
$$
and conjugate the homogeneous systems, $x'=L_tx$ and $\tilde x '=\tilde L_t \tilde x$, with $L_t$ and $\tilde L_t$ Hamiltonian, then necessarily $\gs_t$ is symmetric and $\alpha$ is constant.
Indeed, we can compute
$$
\gR_t'\gR_t^{-1} +a(t)\gR_t L_t\gR_t^{-1}=
\begin{bmatrix}
aC^T-a\alpha^{-1}B \gs_t&\alpha^{-1}aB\\ * & \alpha'\alpha^{-1} I+a\alpha^{-1}\gs_tB-aC
\end{bmatrix}
$$
(the lower left block is not useful). In order for this matrix to be Hamiltonian, we must have 
$$
aC-a\alpha^{-1}\gs_t^TB=aC-a\alpha^{-1}\gs_tB-\alpha'\alpha^{-1}I
$$
which is equivalent to
$$
a(\gs_t^T-\gs_t)=\alpha'B^{-1}.
$$
In this equality, the left hand side is antisymmetric, while the right hand side is symmetric, so both most be null,
which is precisely implying that $\alpha$ is constant and that $\gs_t$ is symmetric.

Proving the direct implication in Proposition \ref{prop-3co} is more delicate. We temporarily work in a more general setting.
Let $L_t=L(t)$ and $\tilde L_t=\tilde L (t)$ be two smooth curves of endomorphisms of a vector space 
$E$ ($=\Rm^d\times (\Rm^d)^*$ in our case), and let $a$ and $\tilde a$ be non-vanishing functions.
We consider a quotient $\pi:E\lto G$ of $E$ (the projection to the first factor $\Rm^d$ in our case) and a subspace $F$ of $E$
(the subspace $\{0\}\times (\Rm^d)^*$) , we assume that $\pi_{|F}=0$.
 Our goal is to give necessary conditions in order that 
the solutions of the two non-homogeneous equations 
\begin{align*}
 x'(t) & = a(t)L_t x(t)+a(t)b(t)\\
 \tilde x'(t)&=\tilde a(t)\tilde L_t \tilde x(t)+\tilde a(t) b(t)
\end{align*}
starting with the initial condition $x(0)=0=\tilde x(0)$ have the same projection on $G$ (meaning that 
$\pi \circ x=\pi \circ \tilde x)$ on the time interval $[0,T]$ for each curve $b:[0,T]\lto F$ smooth and compactly supported in $]0,T[$. Recall that the curve 
$b$ appearing in both differential equations is the same.

Given a linear map $M$ on $E$, we denote by $[M]:F\lto G$ the block $\pi \circ M_{|F}$ of $M$, it is the upper right block in our case.

\begin{lem}\label{lem-pnh}
If $\pi\circ x=\pi \circ \tilde x$ on $[0,T]$, for each curve $b$ smooth and compactly supported in $]0,T[$, then 
$$
[a(t)M_n(t)]= [\tilde a(t) \tilde M_n(t)]
$$
 for each $n\geq 1$ and each $t\in [0,T]$, where $M_n(t)$ is the sequence of curves of matrices defined by $M_1=I$
and $M_{n+1}=M_n'+aM_nL$ (and similarly for $\tilde M_n$).
In the special case where $a\equiv \tilde a$, this implies that 
$
[L_t^3]=[\tilde L_t^3].
$
\end{lem}

In the formulation of the problem, we have included two functions $a$ and $\tilde a$ to make the expressions more symmetric.
It would however be possible to multiply both $a$ and $\tilde a$ by a same non vanishing function without changing the problem.
This would for example allow to assume that $a$ or $\tilde a$ (but, in general, not both) is 
identically equal to $1$

\proof
If $\pi\circ x =\pi \circ \tilde x$ for each curve $b$ smooth and compactly supported in $]0,T[$, then this also holds for
all $b\in L^1([0,T])$, by a standard approximation argument using  that the maps $b\lmto x$ and $b\lmto \tilde x$ are continuous from $L^1$ to $L^{\infty}$.

A straightforward computation  shows that 
\begin{align*}
x'=&ab+ aLx=ab+M_2x\\
x''=&(ab)'+aM_2b+M_3x\\
x^{(3)}=&(ab)''+(aM_2b)'+aM_3b+M_4 x\\
x^{(n)}=&(ab)^{(n-1)}+(aM_2b)^{(n-2)}+ \cdots + a M_n b+ M_{n+1}x.
\end{align*}
We consider a discontinuous curve $b$ which is null on $[0,s[$, and has a prescribed constant value 
$\underline b
$
on $[s, T]$.
The equation $\pi \circ x''(s)=\pi \circ \tilde x''(s)$ (the derivatives are computed on the right of $s$) reduces to 
$$
a(\pi \circ M_2 )\underline b=\tilde a(\pi \circ \tilde M_2 )\underline b
$$
since $x(s)=\tilde x(s)=0$ and $\pi \circ (ab')(s)=0=\pi \circ (\tilde a b)'(s)$. Since this holds for each 
$\underline b$, we have $a[M_2(t)]=\tilde a[\tilde M_2(t)]$ for each $s$ in $]0,T[$ hence in $[0,T]$.

Assuming this equality, we have that $[(aM_2)'(t)]= [(\tilde a\tilde M_2)'(t)]$ hence  the equation
$\pi \circ x'''(s)=\pi \circ \tilde x'''(s)$ 
simplifies to $[a(s)M_3(s)]=[\tilde a (s)M'_3(s)]$. Similarly by induction the equality 
 $\pi \circ x^{(n)}(s)=\pi \circ \tilde x^{(n)}(s)$  reduces to 
 $[a(s)M_n(s)]=[\tilde a (s)M'_n(s)]$.

Assume now that $\tilde a \equiv a$. We have $M_2=aL$, $M_3=(aL)'+(aL)^2$, 
$$
M_4=M_3'+(aL)'(aL)+(aL)^3= M_3'+((aL)^2)'/2+(aL)^3.
$$
Under the hypothesis of the proposition, we have $[aL]=[  a\tilde L]$ hence
$[(aL)']=[(  a\tilde L)']$.
Then, the equality 
$[M_3]=[\tilde M_3]$ reduces to  $[(aL)^2]=[( a \tilde L)^2]$.
Finally, the equality 
$[M_4']=[\tilde M_4']$ implies $[(aL)^3]=[(a\tilde L)^3]$, hence $[L^3]=[\tilde L^3]$.
\qed

\proof We now prove the direct implication in Proposition \ref{prop-3co}.
The first conclusion in Lemma \ref{lem-pnh} implies that 
 $$\tilde a ^2\tilde B_t=[\tilde a^2\tilde L_t]=[a^2L_t]=a^2B_t.
 $$
This implies that the lower right block $(a/\tilde a)\tilde B^{-1}B $ of $\gR_t$ is actually equal to $(\tilde a/a ) I$.

The first block line of 
$L^2$ is 
$$\pi \circ L^2=((C^T)^2-BA, C^TB-BC)$$
and similarly for $\tilde L^2$.
The equality $aM_3=\tilde a _3 \tilde M_3$ can be written
$a [aL]'+a^3[L^2]= \tilde a [\tilde a\tilde L]'+\tilde a^3[\tilde L^2]$
and
$$
a^3(C^TB-BC)-a(aB)'=\tilde a^3(\tilde C^T\tilde B-\tilde B\tilde C)-\tilde a( \tilde a \tilde B)'.
$$
Using that $\tilde a ^2\tilde B=a^2B$, this is equivalent to
\begin{align*}
&\quad aC^T(a^2B)-(a^2B)aC-\tilde a\tilde C^T (a^2B)+ (a^2B)\tilde a \tilde C=
a(aB)'-\tilde a (\tilde a B)'\\
\Leftrightarrow & \quad
a^2B(\tilde a\tilde C-aC) - \big( a^2B(\tilde a\tilde C-aC)  \big)^T=
a(aB)'-\tilde a (\tilde a B)'.
\end{align*}
In this last equality, the matrix on the left is antisymmetric, while the matrix on the right is symmetric
(recall that $B$ is symmetric). The equality thus implies that each of them is null.
The nullity of the left hand side implies that the matrix $B(\tilde a\tilde C-aC)$ is symmetric, hence (using that $B$ is symmetric) so is the transpose
$(\tilde a \tilde C^T-aC^T)B$, and so is 
$$
B^{-1}(\tilde a \tilde C^T-aC^T)=B^{-1}(\tilde a \tilde C^T-aC^T)B B^{-1},
$$
hence (using that $\tilde a ^2\tilde B=a^2B$) so is the lower left block 
$$\gs_t=(\tilde a(t)\tilde B(t))^{-1}\big( a(t)C^T(t)-\tilde a(t) \tilde C^T(t)\big)
$$
of $\gR_t$. 
The nullity of the right hand side, 
$
a(aB)'=\tilde a (\tilde a \tilde B)'
$
 can be rewritten 
 $$ (a^2B)'-a'aB=(\tilde a^2 \tilde B)'-\tilde a' \tilde a \tilde B
 $$
 and, since $a^2B=\tilde a^2 \tilde B$, this is equivalent to $(a'/a -\tilde a'/\tilde a)B=0$ hence to
 $  a'/a =\tilde a'/\tilde a$, which implies that $\tilde a/a$ is constant.
We have proved the first and the second conclusions of the Proposition.

Let us now prove our last conclusion. At this stage, we know that 
$$
\gR_t = \begin{bmatrix}
	I&0\\ \gs_t & (\tilde a /a)  I
\end{bmatrix}.
$$
The computations made in the proof of the converse imply that 
$\gR_t$ conjugates the system $x'=aL_tx+ab$  to the system $y'=\tilde a G_t y +\tilde a b$,
with 
$$
\tilde a (t) G_t=\gR_t'\gR_t +a(t)\gR_t^{-1} L_t\gR_t^{-1}.
$$
By definition of $\gR_t$, the first block line of this matrix $G_t$ is equal to the first block line of $\tilde L$.
Moreover, as seen in the proof of the converse above, the second conclusion implies that the lower right block of $G$
is minus the transpose of its upper left block, so that 
$$
G=\begin{bmatrix}
\tilde C^T & \tilde B\\ -A_G & -\tilde C 
\end{bmatrix},
$$
and we have to prove that $A_G=\tilde A$.

We  apply Lemma \ref{lem-pnh} to the systems $\tilde x'=\tilde a \tilde L \tilde x+ \tilde a b$ and 
$y'=\tilde a G y +\tilde a b$, that is to the pairs of data $(\tilde L, \tilde a)$ and $(G,\tilde a)$.
We obtain that $[G^3]=[\tilde L^3]$.
Using the expression given above for the first block line of $L^2$, we 
expand this   equality to
$$((\tilde C^T)^2-\tilde B \tilde A)\tilde B-(\tilde C^T\tilde B-\tilde B\tilde C)\tilde C=
((\tilde C^T)^2-\tilde B A_G)\tilde B-(\tilde C^T\tilde B-\tilde B\tilde C)\tilde C,
$$
which simplifies to $\tilde BA_G\tilde B= \tilde B \tilde A \tilde B$, hence to $A_G=\tilde A$.
This implies that $G=\tilde L$, and proves the third conclusion of the proposition.
\qed

\section{Reversible points and reversible orbits}\label{sec-rev}

The goal of this  section and the next one  is to prove  Theorem \ref{thm-g5}, 
which will be deduced from the stronger Theorem \ref{thm-twr} below.
We say that $x\in T^*M$ is a two-way point if $x\not \in \Gamma$, and if $\gS(\varphi(t,x))$ is, up to parametrization,
an orbit of $H$ near $t=0$. We say that $(x,u)\in T^*M \times C^{\infty}(M)$ is a two-way point if $x$ is a two way point for $H+u$.
This holds if and only if there exists $\epsilon>0$ such that the equality
$$
d\gS_{\varphi(t,x,u)}\cdot V_H(\varphi(t,x,u),u)=-\ms(\varphi(t,x,u))V_H(\gS\circ \varphi(t,x,u),u)
$$
holds for each $|t|<\epsilon$. We denote by $\mW(\epsilon)$ the set of points $(x,u)$ which satisfy this property, and by
$\mW=\cup_{\epsilon>0} \mW(\epsilon)$ the whole  set of two-way points. This union is increasing, hence it can be made countable. 
It is clear from the definition that $\mW(\epsilon)$ is closed, hence $\mW$ is an $F_{\sigma}$.

If $\theta(t), |t|<\epsilon$ is an orbit segment  made of two-way points (which is equivalent to $(\theta(0),0)\in \mW(\epsilon)$), then
there exists a decreasing diffeomorphism $\sigma $ form $]-\epsilon,\epsilon[$ into an open interval of $\Rm$, satisfying 
 $\sigma(0)=0$, and an orbit segment 
$\tilde \theta$  such that $\gS\circ \theta = \tilde \theta\circ \sigma$. The function $\sigma$ satisfies $\sigma(0)=0$, and  is determined by the equation
 $ \sigma'(t)=-\ms(\theta(t))$.

\begin{defn}\label{def-rp}
We say that the point $\theta$ is  reversible if it is a two way point, and if moreover 
\begin{enumerate}
	\item 
$d\ms_{\theta}\cdot V_H(\theta)=0$,
\item
the reduced linear symmetry $\gR_{\theta}$ is conformally symplectic of factor $-\ms(\theta)$,
\item
denoting by $\theta(t)$ the orbit of $\theta$, we have 
$$
\partial_{t|t=0} \big(\gR_{\gS(\theta)}\circ\tilde  \Phi_{\sigma(t)}^0\circ\gR_{\theta(t)}\circ \Phi_0^t\big) =0
$$
where $\Phi_0^t$ is the reduced linearized flow along $\theta$ and $\tilde \Phi_0^t$ is the reduced linearized flow along $\tilde \theta$.
\end{enumerate}
 \end{defn}

%
%

We are going to prove :

\begin{thm}\label{thm-twr}
	For generic $u\in C^{\infty}(M)$, all two-way points are reversible.
\end{thm}

That Theorem \ref{thm-twr}  implies Theorem \ref{thm-g5} follows from the next Lemma, since all points
of a round trip orbit which do not belong to $\Gamma$  are  two way points.

\begin{lem}\label{lem-rtr}
	Let $\theta(t)$ be a round trip periodic orbit. Assume that all points of the image of $\theta$ which are outside of $\Gamma$ are reversible points. Then the orbit $\theta$ is reversible.
\end{lem}

\textsc{Proof of Lemma \ref{lem-rtr}.} We assume, without loss of generality, that $\theta(0)\in \Gamma$, and denote by $T>0$ the minimal period of the orbit.
There exists one and only one time $\nu\in ]0,T[$ such that $\theta (\nu)\in \Gamma$. 
The first condition of the definition of reversible points implies that $ \ms\circ \theta$ is constant on $]0,\nu[$ and
on $]\nu,T[$. Moreover, we know that $\ms\circ \theta$ is continuous and takes the value $1$ at $\theta(0)$ and $\theta(\nu)$.
We deduce that $\ms\equiv 1$ on $x$, hence that $\nu=T/2$. The time symmetry $\sigma$ associated to the orbit $\theta$ is 
$\sigma (t)=-t$.

For the sequel, we denote by $\gR_t$ the reduction of $d\gS_{\theta(t)}$. Note that $\gR_t\circ \gR_{-t}=Id$.
The maps $\gR_t, t\in ]0,1/2[$ are conformally symplectic of factor $\ms(\theta(t))$, hence at the limit $\gR_0$ and $\gR_{T/2}$
are conformally symplectic of factor $-1$, \textit{i.e.} antisymplectic.

We want to prove the reversibility relation
$
\gR_0\circ \Phi_{T/2}^T\circ \gR_{T/2}\circ \Phi_{0}^{T/2}=Id
$,
which is equivalent to 
$$
\gR_0\circ \Phi_{-T/2}^0\circ \gR_{T/2}\circ \Phi_{0}^{T/2}=Id.
$$
By continuity, it is enough to prove that the reversibility relation
$$
\gR_{-s}\circ  \Phi_{-t}^{-s}\circ\gR_{t}\circ \Phi_s^t=Id.
$$
holds for $s$ and $t$ in $]0,T/2[$.
We denote 
$$\gL_s^t:=
\Phi_s^0 \circ 
\gR_{-s}\circ  \Phi_{-t}^{-s}\circ\gR_{t}\circ \Phi_0^t
=\Phi_s^0
\circ 
\gR_{-s}\circ  \Phi_{-t}^{-s}\circ\gR_{t}\circ \Phi_s^t\circ \Phi_0^s,
$$
it is an endomorphism of the reduced tangent space at $\theta(0)$, which depends smoothly on $s$ and $t$.
Note that $\gL_s^s=Id$, and the identity we try to prove is equivalent to $\gL_s^t\equiv Id$.
We have 
$$
\gL_s^r=\gL_s^t\circ \gL_t^r
$$
for each $s,t,r$ in $]0,T/2[$, as follows from the computation
\begin{align*}
	\gL_s^r 
	&= \Phi_s^0 \circ \gR_{-s}\circ  \Phi_{-r}^{-s}\circ\gR_{r}\circ \Phi_0^r\\
	&=\Phi_s^0 \circ \gR_{-s}\circ \Phi_{-t}^{-s}\circ \Phi_{-r}^{-t}\circ \gR_{r}\circ \Phi_0^r\\
	&=\Phi_s^0 \circ \gR_{-s}\circ \Phi_{-t}^{-s}\circ \gR_t\circ \gR_{-t}\circ \Phi_{-r}^{-t}\circ \gR_{r}\circ \Phi_0^r\\
	&=\Phi_s^0 \circ \gR_{-s}\circ \Phi_{-t}^{-s}\circ \gR_t\circ \Phi_0^t\circ \Phi_t^0\circ  \gR_{-t}\circ \Phi_{-r}^{-t}\circ \gR_{r}\circ \Phi_0^r.
\end{align*}
The third condition of the definition of reversible points implies that
$$\partial_{r|r=t}\gL_t^r=0$$
and in view of the above relation this implies that 
$$\partial_t \gL_s^t\equiv 0,$$ 
and we deduce that $\gL_s^t=Id$ for each $s$ and $t$ in $]0,T/2[$, which implies the desired reversibility relation.
\qed
 
 The proof of Theorem \ref{thm-twr} is based on the following Proposition :
 
 \begin{prop}\label{prop-Wloc}
 Let $\theta(t) =(Q(t),P(t)), |t|<\epsilon $ be an orbit segment made of two way points of  $H$ (hence $(\theta(0),0)\subset \mW(\epsilon)$).
 Assume moreover that $\theta(0)$ is not reversible.
 Then there exists an open neighborhood $\mW_{loc}(\epsilon)$ of 
 $(\theta(0),0)$ in $\mW(\epsilon)$   whose projection on
 $C^{\infty}(M)$ has empty interior.
\end{prop}

\noindent
\textsc{Proof of Theorem \ref{thm-twr}.}
We assume Proposition \ref{prop-Wloc}.
Let $\mW'(\epsilon)$ be the set of elements of $\mW(\epsilon)$ which are not reversible.
  By  Proposition  \ref{prop-Wloc}, applied to $H+u$, 
each point $(\theta,u)\in \mW'(\epsilon)$ is contained in an open subset $\mW_{loc}(\epsilon)$ of $\mW(\epsilon)$ whose projection  on $C^{\infty}(M)$ has empty interior. Moreover, $\mW_{loc}(\epsilon)$ is locally closed, hence it is an $F_{\sigma}$, so its projection 
is an $F_{\sigma}$ with empty interior.
Since $T^*M\times C^{\infty}(M)$ is a separable metric space, the subset $ \mW'(\epsilon)$ can be covered by countably many 
such neighborhoods, so its projection is contained in a countable union of  $F_{\sigma}$ with empty interior, hence, by the Baire property, in an  $F_{\sigma}$ with empty interior.
Then the projection of the set  $\mW'=\cup _{k}\mW'(1/k)$ of two-way points which are not reversible is contained in an $F_{\sigma}$ with empty interior. 
\qed

\noindent
\textsc{Proof of Proposition \ref{prop-Wloc}.}
The principle of the proof is that a point which is stably a two way point has to be a reversible point.
We assume, without loss of generality, that $H(\theta)=0$.
We work in coordinates $(q_0, q_1, \ldots,q_d)=(q_0, q_*)$ of $M$ near $Q(0)=\tilde Q(0)$ such that 
$$ Q(t)=(Q_0(t),0),$$ 
with $Q'(t)>0$, for $t\in ]-\delta, \delta[$ for some $\delta\in ]0,\epsilon[$. Such coordinates exist because $ Q'(0)\neq 0$.
We could obviously impose that $Q_0(t)=t$, but some expressions will appear more symmetric if we don't.
The first coordinate $Q_0$ of the orbit  maps $]-\delta, \delta[$ diffeomorphically  onto an open interval of $\Rm$, and we denote by $\tau $ 
its inverse, so that 
$$Q\circ \tau (r)=re_0=(r,0).$$
The corresponding coordinates of $T^*M$ are 
$$
(q,p)=(q_0,q_1,\ldots,q_d,p_0, p_1, \ldots p_d)=(q_0, q_*, p_0, p_*)
\in \Rm^{d+1}\times \Rm^{d+1}=\Rm\times \Rm^d\times \Rm \times \Rm^d.
$$
The coordinates $x_*:= (q_*,p_*)$ are local symplectic coordinates on the restricted sections $\{q_0=r, H=c\}$, 
hence on the restricted linear sections, hence on the reduced energy levels. 

Since $\theta(0)$ is a two way point, the orbit $\tilde \theta(t)=(\tilde Q(t), \tilde P(t))$ of the point $\gS(\theta)$ satisfies 
$$
\tilde Q(t)=(\tilde Q_0(t), 0)
$$
and we denote by $\tilde \tau$ the inverse of $\tilde Q_0$, it is a decreasing diffeomorphism from a neighborhood of $0$ to $]-\delta, \delta[$. The time symmetry is $\sigma=\tilde \tau \circ Q_0$.

We denote as usual by $\Phi_s^t$, resp. $\tilde \Phi_s^t$,  the reduced linearized flow along the orbit $\theta$, resp. $\tilde \theta$, expressed  in coordinates $x_*$.
It will actually be useful to parameterize orbits by the coordinate $q_0$, and define, in coherence with section \ref{sec-nhls},  $\Psi_{\rho}^r:= \Phi_{\tau(\rho)}^{\tau(r)}$ and
 $\tilde \Psi_{\rho}^r:=\Phi_{\tilde \tau(\rho)}^{\tilde \tau(r)}$. These are the differentials at the orbit of the local transition maps between the restricted sections $\{H=0, q_0=\rho \}$ and $\{H=0, q_0=r \}$, expressed in coordinates $x_*$.
 The matrices $\Psi_{\rho}^r, \tilde \Psi_{\rho}^r$ are symplectic, hence the matrices $L_r, \tilde L_r$ such that
$$
\partial_r \Psi_s^r=\tau'(r)L_r \Psi_s^r
\quad, \quad 
\partial_r \tilde \Psi_s^r=\tilde \tau'(r)L_r \tilde \Psi_s^r
$$
 are Hamiltonian, meaning that they are of the form
$$
L_r=\begin{bmatrix}
	C_r^T & B_r\\-A_r &-C _r
\end{bmatrix}
\quad , \quad 
\tilde L_r=\begin{bmatrix}
	\tilde C_r^T & \tilde B_r\\-\tilde A_r &-\tilde C _r
\end{bmatrix},
$$
with symmetric blocs $A$ and $B$. We will use later the matrices $A_r, B_r, C_r, \tilde A_r, \tilde B_r, \tilde C_r$ defined by the above equalities.
If the coordinates could be chosen such that $\theta_*$ and  $\tilde \theta_* $ are constant, then we would have  $L_r=J\partial^2_{x_*x_*}H(\theta\circ \tau(r))$ and $\tilde L_r=J\partial^2_{x_*x_*}H(\tilde \theta \circ\tilde \tau(r))$.
However, coordinates satisfying these two conditions do not exist in general, and  the above relation between $L$ and $H$ does not necessarily hold,
although only the blocks $B$ remain equal:

\begin{lem}\label{lem-reduc}
	We have
	$$
	B_r= \partial^2_{p_*p_*}H(\theta\circ \tau(r))\quad , \quad \tilde B_r= \partial^2_{p_*p_*}H(\tilde \theta\circ \tilde \tau(r))
	$$ 
	hence they are  positive definite.
\end{lem}

\begin{lem}\label{lem-R}
The expression in coordinates of the reduced linearized symmetry $\gR_{\theta\circ \tau(r)}$ is 
	$$
	\gR_{r}:=
	\begin{bmatrix}
		I&0\\
		(\tilde \tau'(r)\tilde B(r))^{-1}\big(\tau'(r)C^T(r)-\tilde \tau'(r) \tilde C^T(r)\big) &  (\tau'(r)/\tilde \tau'(r)) \tilde B^{-1}(r) B(r)
	\end{bmatrix}.
	$$
\end{lem}

Note that this expression coincides (unsurprisingly) with the matrix called $\gR_t$ before Proposition \ref{prop-3co}
with $a	=\tau'$, $\tilde a  =\tilde \tau '$.
For $x\in T^*M$ near $\{\theta(0),\tilde \theta(0)\}$, $u\in C^{\infty}$ near $0$, $r\in \Rm$ near $0$, we denote by $\psi(r,x,u)\in \Rm^{2d}$ the $x_*$ coordinate of the local intersection
of the orbit of $x$ for $H+u$ with the section $\{q_0=r\}$ (along $\theta$ or along $\tilde \theta$), and set
$\tilde \psi(r,x,u)=\psi(r,\gS(x),u)$.

\begin{lem}\label{lem-diff} Given a potential $u\in C^{\infty}(M)$ vanishing on $Q(]-\delta, \delta[)$, the directional derivatives
	 $y(s):= \partial_{u}\psi(s,\theta, 0)\cdot u$ and $\tilde y(s):= \partial_{u}\tilde \psi(s,\theta, 0)\cdot u$  satisfy the non-homogeneous differential equations 
	$$
	y'(s)= \tau'(s)L(s)y(s)+\tau'(s)b(s)\quad, \quad
	\tilde y'(s)= \tilde \tau '(s)\tilde L(s)\tilde y(s)+\tilde \tau'(s)b(s)
	$$
	where  $b(s)= (0,-\partial_{q_*}u(s,0))\in \Rm^d \times (\Rm^d)^*$. They take the initial values $y(0)=0=\tilde y(0)$.
\end{lem}

We postpone the proof of these Lemma to the next section and  continue the proof of Proposition \ref{prop-Wloc}.
We will apply Proposition \ref{prop-3co} to the two non-homogeneous linear equations appearing in Lemma \ref{lem-diff}, so 
in the notations of Proposition \ref{prop-3co} $a(t)=\tau'(t)$ and $\tilde a (t)= \tilde \tau'(t)$.
If  $\theta(0)$ is not reversible 
we claim that the conclusions of Proposition \ref{prop-3co} are violated on any subinterval of $[0, \delta[$, provided $\delta>0$
is small enough.

The first possibility is that $d\ms\cdot V_H\neq 0$ at $\theta(0)$, 
or equivalently $(\ms\circ \theta)'(0)\neq 0$.
Observing that 
$$ 
\frac{a(r)}{\tilde a(r)}= \frac{\tau' (r)}{\tilde \tau'(r)}=\frac{\tilde Q'_0 \circ \tilde\tau (r)}{Q'_0\circ  \tau(r)}=
-\ms\circ \theta\circ \tau(r),
$$
we deduce that $(a/\tilde a)'(0)\neq 0$ and then we can assume by possibly taking a smaller $\delta>0$ that 
$(a /\tilde a)'(r)\neq 0$
on $[0, \delta[$, hence $\tilde a/a$ is not constant on any subinterval of $]0, \delta[$.

The second possibility is that $\gR_0=\gR_{\theta(0)}$ is not conformally symplectic of factor $-\ms(\theta(0))$, and then we can assume by taking $\delta$ small enough
that $\gR_r$ is not conformally symplectic of factor $-\ms\circ \theta\circ \tau(r) $ for any $r\in [0,\delta[$.

The third and last possibility is that 
$$
\partial_{t|t=0} \big(\gR_{\gS(\theta(0))}\circ\tilde  \Phi_{\sigma(t)}^0\circ\gR_{\theta(t)}\circ \Phi_0^t\big) \neq 0
$$
which is equivalent to 
$$
\partial_{r|r=0} \big(\gR_0^{-1}\circ \tilde \Psi_r^0\circ\gR_{r}\circ \Psi_0^r\big) \neq 0
$$
and then, taking $\delta$ small enough, we deduce that the 
conjugacy equality 
$$
\gR_{r}\circ \Psi_s^r= \tilde \Psi_s^r\circ \gR_{s}
$$
holds only for  $s=  r$ if both $s$ and $r$ are in $]0, \delta[$.

In all cases, the conclusions of Proposition \ref{prop-3co} are violated on all subintervals of $[0, \delta[$.
Fixing a large integer $k$ and a sequence $0<t_1< t_2 < \cdots <t_k<\delta$ of times, we deduce from 
Proposition \ref{prop-3co} applied on each interval $[t_i, t_{i+1}]$  the existence of smooth curves 
$b_i(t)=(0, -\beta_i(t))$ supported in 
$]t_i, t_{i+1}[$, such that the solutions $y_i$, $\tilde y_i$ of the systems 
$$
y'(t)=\tau'(t)L_ty(t)+\tau'(t)b_i(t)\quad, \quad
\tilde y'(t)=\tilde \tau '(t)\tilde L_t\tilde y(t)+\tilde \tau'(t) b_i(t)
$$
emanating from $y_i(t_i)=0=\tilde y_i(t_i)$ have a different projection  at some time $s_i\in ]t_i, t_{i+1}[$.
We consider smooth potentials $u_i\in C^{\infty}(M)$, null on $Q([-\delta, \delta])$, with the property that,
for each $t\in [-\delta, \delta]$, 
$$\partial_{q_*} u_i(te_0)=\beta_i(t).
$$
In view of Lemma \ref{lem-diff}, 
$$
\partial_u \psi (s,\theta,0)\cdot u_i = y_i(s)
\quad,\quad 
\partial_u \tilde \psi (s,\theta,0)\cdot u_i = \tilde y_i(s).
$$
Finally, there exists one coordinate $q_{j_i}$ such that the  $\Rm$-valued function 
$$
\chi_i(x,u):=q_{j_i}\circ\psi(s_i,x,u)-q_{j_i}\circ \tilde \psi(s_i,x,u)
$$
satisfies
$$
\partial_u \chi_i(\theta,0)\cdot u_i \neq 0.
$$
Let $E$ be the  $k$-dimensional vector space of $C^{\infty}(M)$ generated by the potentials $u_i, 1\leq i \leq k$. We define
$$
\chi:  T^*M_{loc}\times C^{\infty}_{loc}\ni(x,u)\lmto (\chi_1(x,u), \ldots, \chi_k(x,u))\in \Rm^k,
$$
where $T^*M_{loc}$ and  $ C^{\infty}_{loc}$ are  open neighborhoods of $x$ and $0$ in $T^*M$ and $C^{\infty}(M)$.
By construction 
$
\partial_u \chi_i(\theta,0)\cdot u_i \neq 0$ 
and 
$
\partial_u \chi_i(\theta,0)\cdot u_j = 0$ for $j >i$, hence 
the restriction to $E$ of the differential $\partial_u\chi (\theta,0)$ is an isomorphism from $E$ to $\Rm^k$.
If a point $(x,u)\in T^*M_{loc}\times C^{\infty}_{loc}$ belongs to $\mW(\epsilon)$, then $\chi(x,u)=0$.
Let us now consider the modified map
$$
\hat \chi :T^*M_{loc}\times C^{\infty}_{loc}\times E_{loc} \ni (x,u,v)\lmto \chi(x,u+v),
$$
where $E_{loc}$ is a small open neighborhood of $0$ in $E$.
By continuity of $\partial_v \hat \chi$, if the neighborhoods are small enough, then $\partial_v\hat \chi$ is an isomorphism at each point.
Then, for each fixed $u\in C^{\infty}_{loc}$, the map 
$\chi_u:(x,v)\lmto \hat \chi (x,u,v)$ is a submersion on $T^*M_{loc}\times E_{loc}$.
This implies that $\chi_u^{-1}(0)$ is a submanifold of codimension $k$ in the finite dimensional manifold 
$T^*M_{loc}\times E_{loc}$, hence of dimension $2d+2$. If $k=\dim E>2d+2$, then we deduce that the projection of $\chi_u^{-1}(0)$ 
on the second factor has no interior in $E$, and in particular there exist arbitrarily small $v\in E$ which do not belong to this projection.
For such a $v$, we have $\chi(x,u+v)\neq 0$ for each $x\in T^*M_{loc}$  hence the potential $u+v$ does not belong to the projection of 
$\mW_{loc}(\epsilon):= \mW(\epsilon)\cap \big(T^*M_{loc}\times C^{\infty}_{loc}\big)$.
\qed

\section{Some computations in coordinates}
The goal of the present section is to prove Lemma \ref{lem-reduc}, \ref{lem-R}, \ref{lem-diff}, which will end the proof of 
Theorem \ref{thm-twr}, hence of Theorem \ref{thm-g5}.
We work in the setting and with the notations introduced in the proof of Proposition \ref{prop-Wloc}.
We are interested in the local transition maps between sections of the form $\{q_0=r, H=0\}$. These sections are symplectic (for the restriction 
of the ambiant symplectic form) and the coordinates $x_*=(q_*,p_*)$ are symplectic local coordinates on them.
These restricted transition maps depend only on the level set $\{H=0\}$. Since $\partial_{p_0} H(\theta(t))=Q'(t)\neq 0$, this energy level can, locally near $\theta$,  be written as a graph
$$
\{H=0\}=\{p_0=-\kappa(q_0,x_*)\},
$$
this property being the definition of the function $\kappa$. Similarly, near $\tilde \theta= \gS(\theta)$, the energy level is given by the equation $q_0+\tilde \kappa(q_0, x_*)=0$.
The dynamics on the energy surface near $\theta$ is, up to positive reparametrization, the same as the one generated by the Hamiltonian $p_0+\kappa$, which is   given by the equations
$$
q_0'=1
\quad, \quad
x_*'=J\partial_{x_*}  \kappa (q_0, x_*).
$$
The local  transition maps near $\theta$ are thus the flow  of the equations
$$
x_*'=J\partial_{x_*}\kappa (t, x_*),
$$
and the linearized transition maps $\Psi$ solve the linearized equations
$$
\partial_r \Psi_{\rho}^r = J\partial^2_{x_*x_*} \kappa (r,\theta_*\circ \tau(r))\Psi_{\rho}^r,
$$
where $\theta_*(t)=x_*\circ \theta(t)= (0, P_*(t))$ is the $x_*$ component of $\theta(t)$. In earlier notations we have
$$
L_r= J\partial^2_{x_*x_*} \kappa (r, \theta_*\circ \tau(r)).
$$

\noindent
\textsc{Proof of Lemma \ref{lem-reduc}.}
From the equation 
$
H(r,-\kappa(r, x_*), x_*)=0,
$
we deduce that 
$$
\partial_{p_0}H(r,-\kappa(r,x_*), x_*) \cdot\partial_{p_*}\kappa(r,x_*)=\partial_{p_*} H(r, -\kappa(r, x_*), x_*).
$$
We differentiate again this equality with respect to $p_*$ at the point $\theta_*\circ \tau (r)$, using that 
\begin{align*}
&	\partial_{p_*}\kappa (r, \theta_*\circ \tau(r) )=0  \quad, \quad 
\partial_{p_*}H(r, \kappa (r,\theta_*\circ \tau(r )),\theta_*\circ \tau(r))=0,\\
&\partial_{p_0}H(r, \kappa (r,\theta_*\circ \tau(r )),\theta_*\circ \tau(r))= Q_0'\circ \tau (r)=1/\tau'(r),
\end{align*}
and get
$$
\partial^2_{p_*p_*}\kappa(r,\theta_*\circ \tau(r) )
=\tau'(r)\partial^2_{p_*p_*}H(r, \kappa (r,\theta_*\circ \tau(r )),\theta_*\circ \tau(r)).
$$
\qed

\noindent
\textsc{Proof of Lemma \ref{lem-R}.}
The symmetry $\gS$ preserves the sections $\{q_0= r, H=0\}$, and it has the form
$$(r, -\kappa(r,q_*,p_*),q_*,p_*)\lmto (r, -\tilde \kappa(r,q_*,\gs_*(r,q_*,p_*)),q_*,\gs_*(r,q_*,p_*)), 
$$
and then 
$$
\gR_r=
\begin{bmatrix}
	I&0\\ \partial_{q_*}\gs(r, \theta_*\circ \tau(r))& \partial_{p_*}\gs(r, \theta_*\circ \tau(r))
	\end{bmatrix}.
$$
The map $\gs_*$ is determined by the equation
$$
\partial_{p_*}\tilde \kappa(r,q_*, \gs_*(r,q_*,p_*))=\partial_{p_*}\kappa(r,q_*, p_*).
$$
Differentiating with respect to $p_*$ at $x_*=\theta_*\circ \tau(r)$ 
and observing that 
$\gs_*(r,\theta_*\circ \tau (r))= \tilde P_*\circ \tilde \tau(r)$ yields
$$
\partial^2_{p_*p_*}\tilde \kappa(r, \tilde \theta_*\circ \tilde \tau (r))
\partial_{p_*}\gs_*(r,  \theta_*\circ \tau (r))
=\partial^2_{p_*p_*} \kappa(r,  \theta_*\circ  \tau (r))
$$
which can be rewritten in terms of the blocs of $\tilde L_r$ as
$$
\tilde \tau'(r) \tilde B(r) \partial_{p_*}\gs_*(r,  \theta_*\circ \tau (r))
=\tau'(r) B(r).
$$
From this equation, we find the stated expression for the lower right block $ \partial_{p_*}\gs_*(r,\theta_*\circ \tau (r))$ of $\gR_r$.
We can also differentiate the equation determining $\gs_*$ with respect to $q_*$ and get
$$
\partial^2_{q_*p_*}\tilde \kappa(r, \tilde \theta_*\circ \tilde \tau (r))+
\partial^2_{p_*p_*}\tilde \kappa(r, \tilde \theta_*\circ \tilde \tau (r))
\partial_{q_*}\gs_*(r,   \theta_*\circ  \tau (r))=
\partial_{q_*p_*}\kappa (r, \theta_*\circ \tau(r))
$$
which can be rewritten 
$$
\tilde \tau'(r)\tilde C^T(r)+\tilde \tau'(r)\tilde B(r)\partial_{q_*}\gs_*(r,  \theta_*\circ \tau (r))=\tau'(r)C^T(r),
$$
from which we obtain the lower left block.
\qed

\noindent
\textsc{Proof of Lemma \ref{lem-diff}.}
We fix the potential $u$ and  denote by $\kappa(q_0, x_*,\epsilon)$ the function such that the energy level $\{H+\epsilon u=0\}$ is locally the graph 
$p_0=-\kappa$. 
Then $\psi$ solves the equation
$$
\partial_s \psi (s,x, \epsilon u)= J\partial_{x_*}\kappa ( s, \psi(s,x,\epsilon u),\epsilon u).
$$
Taking the derivative with respect to $\epsilon$ at $x=\theta, \epsilon=0$  and recalling that $\psi(s,\theta,0)=\theta_*\circ \tau(s)$, we obtain 
$$
y'(s)=J\partial^2_{x_*x_*}\kappa (s,\theta_*\circ \tau(s),0)y(s)+J\partial^2_{\epsilon x_*}\kappa(s,\theta_*\circ \tau(s),0).
$$
We have 
$J\partial^2_{x_*x_*}\kappa (s,\theta_*\circ \tau(s),0)
=\tau'(s)L(s)$ 
by definition, and we now have to prove that 
$$J\partial^2_{\epsilon x_*}\kappa(s,\theta_*\circ \tau(s),0)=\tau'(s) b(s).
$$
From the equation 
$$
H(s,-\kappa(s, x_*, \epsilon), x_*)+\epsilon u(s,q_*)=0,
$$
we deduce that
$$
\partial_{p_0}H(s,-\kappa( s, x_*,\epsilon),x_*  ) \partial_{x_*} \kappa ( s, x_*, \epsilon)=
\epsilon \partial_{x_*}u(s,q_*)
+
\partial_{x_*}H (s,-\kappa( s, x_*, \epsilon),x_*  ).
$$
Here  we consider $u$ both as a function of $q$ and as a function of $x$ which depends only on $q$,
so that $\partial_{x_*}u=(\partial_{q_*}u, 0)$. 
Since $u(s,0)=0$, we observe that the equation 
$$
H(s, -\kappa, q_*, p_*)+\epsilon u(s,q_*)=0
$$ 
which defines $\kappa(s,q_*, p_*, \epsilon)$ does not actually depend on $\epsilon$ provided $q_*=0$, and this implies that 
$\kappa(s,0, p_*,\epsilon)$ does not depend on $\epsilon$, and in particular 
$\kappa(s, \theta_*\circ \tau(s), \epsilon)$ does not depend on $\epsilon$.
With this observation in mind, we 
differentiate the above equality  with respect to $\epsilon$ at $\epsilon=0$, $x_*=\theta_*\circ \tau(s)$  and get 
$$
\partial_{p_0}H(\theta \circ \tau(s))\partial^2_{\epsilon x_*} \kappa ( s,  \theta_*\circ \tau(s),0)= \partial_{x_*}u(s,0)
$$
which can be rewritten
$$
\partial^2_{\epsilon x_*} \kappa (s, \theta_*\circ \tau(s),0)= \tau'(s) (\beta(s),0).
$$
\qed

\appendix

\section{A soft implicit function Lemma}\label{ap}

In this text, the perturbation parameter $u$ belongs to a Fréchet space, and we want to avoid using differential calculus in Fréchet spaces. So we use the following elementary version of the implicit function theorem, where $U$ is a metric space. The proof is easy and left to the reader.

\begin{lem}
	Let $f(x,u):\Rm^n\times U\lto \Rm^n$ be a function which is Fréchet differentiable in $x$ for each $(x,u)$, and such that 
the Fréchet differential $\partial_xf:\Rm^n\times U\lto \mL(\Rm^n,\Rm^n)$ is jointly continuous. Assume that $f(0,u_0)=0$, where $u_0$ is some point of $U$, and that $\partial_x f(0,u_0)$ is invertible. Then there exists open neighborhoods $\Rm^n_{loc}$ and $U_{loc}$ of $0$ and $u_0$ and a continuous function $X: U_{loc}\lto \Rm^n_{loc}$ such that, for each $u\in U_{loc}$,
$X(u)$ is the only solution of the equation $f(.,u)=0$ contained in $\Rm^n_{loc}$.
\end{lem}


\begin{thebibliography}{FP}
	
	\bibitem{A63} \textsc{R. Abraham}, Transversality in manifolds of mappings, Bull. A.M.S. \textbf{69}, vol 4 (1963), 470-474.
	
\bibitem{A68}	\textsc{R. Abraham}, Bumpy metrics, Global analysis (Berkeley, Calif., 1968), Proc. Sympos. Pure Math., \textbf{14}, (1970), 1-3.

	
\bibitem{A82} \textsc{D. V. Anosov}, Generic properties of closed geodesics. (Russian) Izv. Akad. Nauk SSSR Ser. Mat. \textbf{46} (1982), no. 4, 675–709, 896; translation : Mathematics of the USSR-Izvestiya, \textbf{21} (1983), No. 1


\bibitem{AB}
\textsc{S. Aslani, P. Bernard,}
Normal Form Near Orbit Segments of Convex Hamiltonian
Systems, IMRN (2021) doi:10.1093/imrn/rnaa344.

\bibitem{AB2}
\textsc{S. Aslani, P. Bernard,} Bumpy metric theorem in the sense of Ma\~{n}e for non-convex Hamiltonians, preprint.


	\bibitem{BK}
\textsc{S. Bolotin, V. V. Kozlov}, Libration in systems with many degrees of freedom, PPM \textbf{42} (1978) no. 2, 245-250.








\bibitem{Co10} \textsc{G. Contreras}, Geodesic flows with positive topological entropy, twist maps and hyperbolicity, Ann. of Math. \textbf{172} (2010), 761-808.


\bibitem{D} \textsc{R. L. Devaney}, Reversible Diffeomorphisms and Flows,  Trans. A.M.S. 
\textbf{218} (1976),  89-113.  

\bibitem{dG} \textsc{M. de Gausson}, Symplectic geometry and Quantum Mechanics, Springer (2006).



	\bibitem{K76}
\textsc{V. V. Kozlov}, The principle of least action and periodic solutions in problems of classical mechanics, PPM \textbf{40} (1976) no. 3, 399-407.

\bibitem{K63}
\textsc{I. Kupka} Contribution à la théorie des champs génériques. (French) Contributions to Differential Equations \textbf{2} (1963), 457–484.

\bibitem{LRR}
\textsc{A. Lazrag, L. Rifford, R. Ruggiero}, Franks' lemma for C2-Mañé perturbations of Riemannian metrics and applications to persistence. J. Mod. Dyn. \textbf{10} (2016), 379–411. 


\bibitem{O08}
\textsc{E. Oliveira}, Generic properties of Lagrangians on surfaces: the Kupka-Smale theorem. DCDS \textbf{21}, no. 2 (2008), 551–69.

\bibitem{RR}
\textsc{L. Rifford, R. Ruggiero}, Generic properties of closed orbits of Hamiltonian flows from Ma\~né's viewpoint. Int. Math. Res. Not., \textbf{22} (2012), 5246–5265. 

\bibitem{CR1}
\textsc{C. Robinson}, Generic properties of conservative systems I. Am. Journ. Maths. \textbf{92} (1970), 562-603.

\bibitem{CR2}
\textsc{C. Robinson}, Generic properties of conservative systems II. Am. Journ. Maths. \textbf{92} (1970), 897-906.

\bibitem{S63}
\textsc{S. Smale},
Stable manifolds for differential equations and diffeomorphisms.
Ann. Scuola Norm. Sup. Pisa Cl. Sci. (3) \textbf{17} (1963), 97–116.







\end{thebibliography}
\end{document}